\let \Bbb \bf
\magnification \magstephalf
%\vsize 9.6truein
\overfullrule0pt

\def\ifundefined#1{\expandafter\ifx\csname#1\endcsname\relax}

\newcount\sectionnumber
\newcount\equationnumber
\newcount\thnumber

\nopagenumbers
\headline={\ifnum\pageno=1 \hfill \else\hss{\tenrm--\folio--}\hss \fi}

\def\assignnumber#1#2{%
  \ifundefined{#1}\relax\else\message{#1 already defined}\fi
  \expandafter\xdef\csname#1\endcsname
  {\if-\the\sectionnumber\else\the\sectionnumber.\fi\the#2}}%
\def\secassignnumber#1#2{%
  \ifundefined{#1}\relax\else\message{#1 already defined}\fi
  \expandafter\xdef\csname#1\endcsname{\the#2}}%
%
% macros on section numbers
%
\def\secname#1{\relax
  \global\advance\sectionnumber by 1
  \secassignnumber{S#1}\sectionnumber
  \csname S#1\endcsname}
\def\Sec#1 #2 {\vskip0pt plus.1\vsize\penalty-250\vskip0pt plus-.1\vsize
  \bigbreak\bigskip
  \equationnumber0\thnumber0
  \noindent{\bf \secname{#1}. #2}\par
  \nobreak\smallskip\noindent}

\def\sectagelse#1{\ifundefined{S#1}\message{S#1 undefined}{\sl #1}%
  \else\csname S#1\endcsname\fi}
%
% macros on equation numbers
%
\def\eq#1{\relax
  \global\advance\equationnumber by 1
  \assignnumber{EN#1}\equationnumber
  {\rm (\csname EN#1\endcsname)}}

\def\eqtag#1{\ifundefined{EN#1}\message{EN#1 undefined}{\sl (#1)}%
  \else(\csname EN#1\endcsname)\fi}
%
% macros on theorem numbers
%
\def\thname#1{\relax
  \global\advance\thnumber by 1
  \assignnumber{TH#1}\thnumber
  \csname TH#1\endcsname}

\def\Assumption#1 {\bLP{\bf Assumption \thname{#1}.}\quad}
\def\Cor#1 {\bLP{\bf Corollary \thname{#1}}\quad}
\def\Def#1 {\bLP{\bf Definition \thname{#1}}\quad}
\def\Example#1 {\bLP{\bf Example \thname{#1}}\quad}
\def\Lemma#1 {\bLP{\bf Lemma \thname{#1}}\quad}
\def\Prop#1 {\bLP{\bf Proposition \thname{#1}}\quad}
\def\Remark#1 {\bLP{\bf Remark \thname{#1}}\quad}
\def\Theor#1 {\bLP{\bf Theorem \thname{#1}}\quad}

\def\thtag#1{\ifundefined{TH#1}\message{TH#1 undefined}{\sl #1}%
  \else\csname TH#1\endcsname\fi}

\def\Proof{\sLP{\bf Proof}\quad}

%
% other macros
%
\hyphenation{pre-print}

\let\bPP=\bigbreak
\def\br={\hfil\break}
\def\LP{\par\noindent}
\def\sLP{\smallbreak\noindent}
\def\mLP{\medbreak\noindent}
\def\bLP{\bigbreak\noindent}

\def\wrt{with respect to }

\def\halmos{\hbox{\vrule height0.31cm width0.01cm\vbox{\hrule height
 0.01cm width0.3cm \vskip0.29cm \hrule height 0.01cm width0.3cm}\vrule
 height0.31cm width 0.01cm}}
\def\hhalmos{{\unskip\nobreak\hfil\penalty50
	\quad\vadjust{}\nobreak\hfil\halmos
	\parfillskip=0pt\finalhyphendemerits=0\par}}
\def\:{:\allowbreak }
\def\refitem{\smallbreak\item}

\def\ga{\gamma}
\def\de{\delta}
\def\ep{\varepsilon}

\def\la{\lambda}
\def\si{\sigma}

\def\De{\Delta}

\def\CC{{\Bbb C}}

\def\RR{{\Bbb R}}

\def\ZZ{{\Bbb Z}}

\def\FSA{{\cal A}}
\def\FSB{{\cal B}}

\def\FSF{{\cal F}}

\def\FSO{{\cal O}}

\def\iy{\infty}
\def\id{{\rm id}}

\def\Zplus{\ZZ_{\,+}}

\def\LHS{left-hand side}
\def\RHS{right-hand side}
\def\half{{\scriptstyle{1\over2}}}

\def\rep{representation}
\let\pa=\partial

\let\ten=\otimes
\let\wt=\widetilde

\def\qchoose#1#2#3{\left[#1\atop#2\right]_#3}
\let\Imply\Longrightarrow
\def\Cqxy{\CC_q[x,y]}
\def\Cqixy{\CC_q[[x,y]]}
\def\Cqhxy{\CC_{qHeis}[[x,y,c]]}

\def\mediumint{{\textstyle\int}}

\font\titlefont=cmr10 at 17.28truept
\font\authorfont=cmr10 at 14truept

\centerline
{\titlefont
Special functions and q-commuting variables}
\bPP
\centerline{\authorfont Tom H. Koornwinder}
\bPP
\centerline{Department of Mathematics}
\centerline{University of Amsterdam}
\centerline{Plantage Muidergracht 24}
\centerline{1018 TV Amsterdam, The Netherlands}
\centerline{email {\tt thk@wins.uva.nl}}
\vskip 1truecm
\noindent
{\sl Abstract:}\quad
This paper is mostly a survey, with a few new results.
The first part deals with functional equations for q-exponentials,
q-binomials and q-logarithms in q-commuting variables and more generally
under q-Heisenberg relations. The second part discusses translation invariance
of Jackson integrals, q-Fourier transforms and the braided line.

\bLP
{\sl Last modified:}\quad August 26, 1996

\bLP
{\sl Note:}\quad
Report No.\ 1, Institut Mittag-Leffler, Djursholm, Sweden, 1995/96;
q-alg/9608008;
\LP
to appear in ``Special Functions, q-Series and Related Topics'',
The Fields Institute Communications Series.

\Sec{1} {Introduction}
Identities for special functions often involve several variables, even if
the special function itself depends on only one variable.
In general these variables are real or complex, so they commute with each other.
The theory of quantum groups has been quite successful in producing identities
for $q$-special funcions, in particular addition formulas,
see e.g.\ the survey in Koelink [21, Section 1] and further references given
there.
Although quantum groups themselves abound of
non-commuting variables satisfying certain relations, one usually does not
find back a similar type of variables in the resulting $q$-special function
identities.
I want to advertize here that special function identities in
non-commuting variables should be studied more extensively and
systematically. They often provide more elegant formulas than the corresponding
identities in commuting variables, and they may be closer to a quantum
group theoretical origin and therefore have more canonical properties.
Another feature (which may be evaluated in a positive or negative sense)
is that such identities are often more algebraic and formal in spirit
and further away from Weierstrass type analysis
than the identities in commuting variables.
The most interesting and challenging cases with non-commuting variables occur
when formal infinite series and convergent infinite series mix with each other.
One has to be extremely careful there in order to avoid paradoxes, see
Section 9.

The present paper surveys (and extends a little) special function theory
involving
$q$-commuting variables $x$ and $y$ (i.e., satisfying the relation
$xy=qyx$ with $q$ complex, usually taken between 0 and 1).

The contents are as follows.
In Section 2 we discuss Sch\"utzenberger's $q$-binomial formula.
Sections 3  deals with various functional equations for $q$-exponentials,
and Section 4 gives some extensions of these results to $q$-Heisenberg cases.
Section 5 describes possible equivalence with formulas involving commuting
variables, via the operational interpretation.
In Section 6 we discuss the $q$-logarithm.
The next four sections are much inspired by the paper by Kempf \& Majid
[19]. We discuss translation invariance under a
$q$-commuting translation variable for Jackson integrals over a finite interval
(Section 7) and over the interval $(-\iy,\iy)$ (Section 9).
In Section 8 we introduce a
$q$-Fourier transform pair in connection with discrete $q$-Hermite polynomials.
While this is in commuting variables, it is related to two types of
$q$-Fourier transforms involving non-commuting variables which have been
studied, respectively, by Kempf \& Majid [19] and
Finkelstein \& Marcus [12].
A deeper understanding of many of these results can be obtained by means
of Majid's [27]
braided quantum groups, in particular the braided line.
This is the topic of Section 10. Finally, two further directions are very
briefly indicated in Section 11.

\mLP
{\sl Conventions}\quad
$\Zplus$ will denote the set of nonnegative integers.
>From Section 3 on,
whenever we work with $q$, it is supposed that $0<q<1$, unless said otherwise.
The notation for $q$-hypergeometric series follows the book by Gasper \&
Rahman [6].

\mLP
\Sec{2} {The q-binomial formula}
Newton's binomial formula says:
$$
(x+y)^n=\sum_{k=0}^n{n\choose k}y^{n-k}x^k,\quad n\in\Zplus.
\eqno\eq{1}
$$
Here it is implicitly understood that $x$ and $y$ commute: $xy=yx$.
A $q$-analogue of \eqtag{1} for {\sl q-commuting variables} $x,y$, i.e.,
satisfying the relation
$$
xy=qyx
\eqno\eq{2}
$$
for some $q\in\CC$,
first appeared
in literature in Sch\"utzenberger [32],
see also Cigler [7, (7)]:

\Prop{1} ({\sl q-binomial formula})\quad
Let $q\in\CC$.
Let $\Cqxy$ be the complex associative algebra with 1 generated by
$x$ and $y$ and with relation \eqtag{2}.
Then the following identity is valid in the algebra $\Cqxy$:
$$
(x+y)^n=\sum_{k=0}^n \qchoose nkq y^{n-k} x^k,\quad
n\in\Zplus.
\eqno\eq{3}
$$
\bPP
Here we used the {\sl $q$-binomial coefficient}
$$
\qchoose nkq:=
{(q;q)_n\over(q;q)_k\,(q;q)_{n-k}}=
{(1-q^n)(1-q^{n-1})\ldots(1-q^{n-k+1})\over
(1-q)(1-q^2)\ldots(1-q^k)}\,,
\eqno\eq{5}
$$
while
the {\sl $q$-shifted factorial} is given by
$$
(a;q)_k:=(1-a)(1-qa)\ldots(1-q^{k-1}a),\quad a\in\CC,\;k\in\Zplus.
\eqno\eq{8}
$$
The recurrence relations below show that the
$q$-binomial coefficient \eqtag{5} is a polynomial in $q$ and therefore
remains meaningful for $q$ being a root of unity.

We will give a constructive proof of Proposition \thtag{1} which goes
essentially back to
Polya \& Alexanderson [29]
and which was later written down by Askey [2].
It is straightforward that the monomials $y^lx^k$ ($k,l\in\Zplus$)
form a basis for $\Cqxy$ considered as a linear space and that
$(x+y)^n$ will have a unique expansion of the form
$$
(x+y)^n=\sum_{k=0}^n c_{n,k}\,y^{n-k}x^k
\eqno\eq{4}
$$
with the coefficients $c_{n,k}$ also depending on $q$.
It follows immediately from \eqtag{4} that $c_{n,0}=1=c_{n,n}$.
Also, expansion of $(x+y)^n=(x+y)^{n-1}(x+y)$ and
$(x+y)^n=(x+y)(x+y)^{n-1}$, respectively, yields for $n>k>0$:
$$
c_{n,k}=q^kc_{n-1,k}+c_{n-1,k-1},\qquad
c_{n,k}=c_{n-1,k}+q^{n-k}c_{n-1,k-1}.
\eqno\eq{6}
$$
Elimination of $c_{n-1,k}$ from these two recurrence equations leaves
us with the two-term recurrence
$$
c_{n,k}={1-q^n\over 1-q^k}\,c_{n-1,k-1}.
$$
Iteration of this last recurrence yields the \RHS\ of \eqtag{5}.
This proves \eqtag{3}.

The advantage of this proof is that it is constructive.
If one just wants to prove by induction with respect to $n$
that \eqtag{4} holds with $c_{n,k}$ being given
by \eqtag{5} then it is sufficient to have only one of the recurrences in
\eqtag{6}. This is the way in which one usually works in case $q=1$,
where the $q$-binomial formula \eqtag{3} by continuity becomes the
binomial formula \eqtag{1}. In that case the two recurrences in \eqtag{6}
coincide and it is not possible to get a two-term
recurrence formula by elimination.

A second observation, due to Andrews, and written down in Askey [3],
is that
the $q$-binomial formula \eqtag{3} is equivalent to an identity in
commuting variables. Note that, if the generators $x,y$ of $\Cqxy$
satisfy the $q$-commutation relations \eqtag{2}, then $-yx,y$ also
satisfy these relations:
$$
(-yx)y=qy(-yx).
$$
Hence, we get from \eqtag{3} that
$$
(-yx+y)^n=\sum_{k=0}^n
\qchoose nkq y^{n-k}(-yx)^k.
\eqno\eq{7}
$$
The \LHS\ of \eqtag{7} equals
$$
\eqalignno{
&y(1-x)y(1-x)\ldots y(1-x)\quad\hbox{($2n$ factors)}
\cr
&\quad=y^n(1-q^{n-1}x)\ldots(1-qx)(1-x)=y^n(x;q)_n.
\cr}
$$
(Note that the definition \eqtag{8} of $q$-shifted factorial remains
valid for $a$ in any complex associative algebra with 1.)$\;$
As for the \RHS\ of \eqtag{7} note that
$$
y^{n-k}(-yx)^k=(-1)^kq^{\half k(k-1)}y^nx^k
$$
and
$$
\qchoose nkq=(-1)^k\,q^{-\half k(k-1)}\,q^{nk}\,{(q^{-n};q)_k\over(q;q)_k}\,.
\eqno\eq{45}
$$
Hence the identity \eqtag{7} can be equivalently written as
$$
y^n(x;q)_n=y^n\sum_{k=0}^n{(q^{-n};q)_k\over(q;q)_k}(q^nx)^k.
$$
Because of the earlier observation about a basis of monomials for $\Cqxy$
we conclude that
$$
(x;q)_n=\sum_{k=0}^n{(q^{-n};q)_k\over(q;q)_k}(q^nx)^k.
$$
This is still an identity in the algebra $\Cqxy$.
However, we can map it to an identity in $\CC$ by using the algebra
homomorphism $\pi\colon\Cqxy\to\CC$ such that, for some $z\in\CC$,
$\pi(x)=q^nz$ and $\pi(y)=0$. This yields the formula giving
the evaluation of a terminating $q$-binomial series:
$$
(q^{-n}z;q)_n=\sum_{k=0}^n{(q^{-n};q)_k\over(q;q)_k}\,z^k,\quad
n\in\Zplus.
\eqno\eq{9}
$$
This formula is well-known, see [16, (II.4)].
It was earlier observed by Cigler [7, (8)] that formula \eqtag{9}
follows from \eqtag{3}. He used an operational interpretation of \eqtag{3}.
We will discuss such interpretations in a later section.
Then we will also see that \eqtag{9} is in fact equivalent to \eqtag{3}.

\Sec{3} {Identities for $q$-exponentials with q-commuting variables}
The two {\sl $q$-exponentials} are defined by
$$
\eqalignno{
e_q(z)&:=\sum_{n=0}^\iy {z^n\over(q;q)_n}={1\over(z;q)_\iy}\,,&\eq{10}
\cr
E_q(z)&:=\sum_{n=0}^\iy {q^{\half n(n-1)}z^n\over(q;q)_n}=
(-z;q)_\iy\,.&\eq{11}
\cr}
$$
Here we assume $0<q<1$ and $(a;q)_\iy$ is defined as the
(convergent) infinite product:
$$
(a;q)_\iy:=\prod_{j=0}^\iy(1-q^j a).
$$
For convergence of the infinite series in \eqtag{10}, \eqtag{11} with
$z\in\CC$ we need $|z|<1$ in \eqtag{10}. However, because of its product \rep,
$e_q$ has an analytic continuation to $\CC\backslash\{q^{-k}\mid k\in\ZZ\}$.
See
[16, Section 1.3]
for the proofs of the second equalities in \eqtag{10} and \eqtag{11}.
The two $q$-exponential series can also be considered as formal power series
in the formal variable $z$. Of course, no convergence condition is needed
in that case.
>From \eqtag{10}, \eqtag{11} we have
$$
e_q(z)\,E_q(-z)=1
\eqno\eq{14}
$$
for $|z|<1$.
>From the second equalities in \eqtag{10} and \eqtag{11} we see that
$$
e_q(qz)=(1-z)\,e_q(z),\quad
E_q(z)=(1+z)\,E_q(qz).
\eqno\eq{13}
$$
In general, algebraic identities for convergent power series
remain valid for the corresponding formal power series.
In particular, this applies to \eqtag{14} and \eqtag{13}.

Fix $q\in(0,1)$ and let $\Cqixy$ be the complex associative algebra with 1
of formal power series
$$
\sum_{k,l=0}^\iy c_{k,l}\,y^lx^k
$$ 
with arbitrary complex coefficients $c_{k,l}$ and
where $x,y$ satisfy relation \eqtag{2}, i.e. $xy=qyx$.
The following Proposition generalizing the classical functional equation
$e^xe^y=e^{x+y}$ for commuting variables $x,y$, was given first
by Sch\"utzenberger [32], see also Cigler [7, (10)].

\Prop{2}
In the algebra $\Cqixy$ we have the identities
$$
\eqalignno{
e_q(x+y)&=e_q(y)e_q(x),&\eq{12}
\cr
E_q(x+y)&=E_q(x)\,E_q(y).&\eq{38}
\cr}
$$
\Proof
By means of the $q$-binomial formula \eqtag{3} we can write the \LHS\
of \eqtag{12} as an element of $\Cqixy$, and next rewrite it as the \RHS:
$$
\eqalignno{
&e_q(x+y)=\sum_{n=0}^\iy\sum_{k=0}^n{1\over(q;q)_n}\,\qchoose nkq y^{n-k}x^k
=\sum_{n=0}^\iy\sum_{k=0}^n {1\over (q;q)_{n-k}(q;q)_k}\,y^{n-k}x^k
\cr
&=\sum_{k,l=0}^\iy {1\over (q;q)_l(q;q)_k}\,y^lx^k
=e_q(y)e_q(x).
\cr}
$$
This settles \eqtag{12}. For the proof of \eqtag{38} note that,
in view of \eqtag{14}, $E_q(x+y)$ is a left and right inverse to
$e_q(-x-y)$, and $E_q(x)E_q(y)$ is a left and right inverse to
$e_q(-y)e_q(-x)$. Now use \eqtag{12}.\hhalmos

\bPP
The reader is warned that the apparent symmetry in $x$ and $y$ in the left
hand side of \eqtag{12} does not allow to conclude that
$e_q(x+y)=e_q(x)e_q(y)$, since the relation $xy=qyx$ is not symmetric
in $x$ and $y$.
The next Proposition gives a formula for $e_q(x)e_q(y)$ in the
algebra $\Cqixy$, i.e.\ for the \RHS\ of \eqtag{12} with the order
of the two factors interchanged. It is a special case of a more general
result given first in operational form by Rogers [30],
which we will discuss in the next section.
See also Gelfand \& Fairlie [17, (46)],
Floreanini \& Vinet [13, formula (23d)],
Faddeev \& Volkov
[11, p.314, formula (2)],
Faddeev \& Kashaev [10, Section 2],
A. N. Kirillov [20, Section 2.5, Lemma 9],
and McDermott \& Solomon [28, (10)].

\Prop{3}
In the algebra $\Cqixy$ we have the identities
$$
\eqalignno{
e_q(x)\,e_q(y)
&=e_q(y-yx)\,e_q(x)&\eq{15}
\cr
&=e_q(x+y-yx)&\eq{18}
\cr
&=e_q(y)\,e_q(-yx)\,e_q(x)&\eq{19}
\cr
&=e_q(y)\,e_q(x-yx),&\eq{20}
\cr
E_q(y)\,E_q(x)&=E_q(x+y+yx),&\eq{39}
\cr
&=E_q(x)\,E_q(yx)\,E_q(y).&\eq{115}
\cr}
$$
\Proof
Because $e_q(x)$ is invertible as a formal power series
(cf.\ \eqtag{14}), formula \eqtag{15} can equivalently be stated as
$$
e_q(x)\,e_q(y)\,e_q(x)^{-1}=e_q(y-yx).
\eqno\eq{16}
$$
For any two formal power series $f(z)$ and $g(z)$,
with $f(z)$ being invertible and $g(z)=\sum_{k=0}^\iy c_kz^k$, we have
$$
f(x)g(y)f(x)^{-1}=
\sum_{k=0}^\iy c_k\,f(x)\,y^k\,f(x)^{-1}=
\sum_{k=0}^\iy c_k\,\bigl(f(x)yf(x)^{-1}\bigr)^k=
g\bigl(f(x)\,y\,f(x)^{-1}\bigr)
$$
as identities in $\Cqixy$.
In particular,
$$
e_q(x)\,e_q(y)\,e_q(x)^{-1}=e_q\bigl(e_q(x)\,y\,e_q(x)^{-1}\bigr).
\eqno\eq{17}
$$
Now, since $xy=qyx$ and by \eqtag{13} we have
$$
e_q(x)\,y\,e_q(x)^{-1}=
y\,e_q(qx)\,e_q(x)^{-1}=
y\,(1-x)\,e_q(x)\,e_q(x)^{-1}=
y\,(1-x).
\eqno\eq{113}
$$
Together with \eqtag{17} this settles \eqtag{16} and hence \eqtag{15}.

Next it follows from Proposition \thtag{2} that the \RHS\ of \eqtag{15}
equals \eqtag{18}, since
$x(y-yx)=q(y-yx)x$.
The equalities \eqtag{19} and \eqtag{20} also follow by
application of Proposition \thtag{2}.
Finally, \eqtag{39} and \eqtag{115} follow from \eqtag{18} and \eqtag{19}
in a similar way as we obtained \eqtag{38}.
\hhalmos

\bPP
Note that the equalities \eqtag{15}--\eqtag{39} reduce to the
classical identities $e^xe^y=e^ye^x=e^{x+y}$ if we replace
$x,y$ by $(1-q)x,(1-q)y$ and let $q\uparrow1$
on using
$$
\lim_{q\uparrow 1} e_q((1-q)z)=e^z=\lim_{q\uparrow 1} E_q((1-q)z),
\eqno\eq{90}
$$
cf.\ [16, (1.3.17)].

\bPP
As a corollary to both Proposition \thtag{2} and Proposition \thtag{3}
we have a functional equation  for the {\sl q-Gaussians}
$$
g_q(x):=e_{q^2}(-x^2),\quad
G_q(x):=E_{q^2}(-x^2)
\eqno\eq{66}
$$
in $q$-commuting variables.
First note that that, for $z\in\CC$ with $|z|<1$, we have
$$
e_{q^2}(-z^2)={1\over(-z^2;q^2)_\iy}
={1\over(iz;q)_\iy(-iz;q)_\iy}=e_q(iz)\,e_q(-iz).
$$
Thus the equality
$$
e_{q^2}(-z^2)=e_q(iz)\,e_q(-iz)
\eqno\eq{67}
$$
is also valid for arbitrary real $z$ or as an identity for formal power
series.

\Cor{14}
In the algebra $\Cqixy$ we have the identities
$$
\eqalignno{
e_{q^2}(-(x+y)^2)&=e_{q^2}(-y^2)\,e_q(-yx)\,e_{q^2}(-x^2),&\eq{68}
\cr
E_{q^2}(-(x+y)^2)&=E_{q^2}(-x^2)\,E_q(-yx)\,E_{q^2}(-y^2).&\eq{133}
\cr}
$$

\Proof
We apply first \eqtag{67}, next Proposition \thtag{2} (twice),
next Proposition \thtag{3} and finally \eqtag{67} again (twice):
$$
\eqalignno{
&e_{q^2}(-(x+y)^2)
=e_q(i(x+y))\,e_q(-i(x+y))
=e_q(iy)\,e_q(ix)\,e_q(-iy)\,e_q(-ix)
\cr
&=e_q(iy)\,e_q(-iy)\,e_q(-yx)\,e_q(ix)\,e_q(-ix)
=e_{q^2}(-y^2)\,e_q(-yx)\,e_{q^2}(-x^2).
\cr}
$$
This yields \eqtag{68}. Then formula \eqtag{133} follows by taking the inverse
on both sides and replacing $x,y$ by $ix,iy$, respectively.\hhalmos

\bPP
Let us next discuss generalizations of the previous results in this section
for the case of non-terminating $q$-binomial series
$$
{}_1\phi_0(a;;q,z):=
\sum_{k=0}^\iy{(a;q)_k\over(q;q)_k}\,z^k=
{(az;q)_\iy\over (z;q)_\iy},\quad
|z|<1,\;a\in\CC,
\eqno\eq{33}
$$
see [16, (II.3)].
Formula \eqtag{33} can be rewritten as
$$
{}_1\phi_0(a;;q,z)=E_q(-az)\,e_q(z),
\eqno\eq{37}
$$
and this remains valid as an identity for formal power series.
The termwise limit for $q\uparrow 1$ of
${}_1\phi_0(q^a;;q,z)$ is
$$
{}_1F_0(a;;z):=\sum_{k=0}^\iy {(a)_k\over k!}\,z^k=(1-z)^{-a}.
$$
The functional equation
$(1-x)^{-a}(1-y)^{-a}=(1-x-y+xy)^{-a}$
in commuting variables $x,y$
can equivalently be written as
$$
{}_1F_0(a;;x)\,{}_1F_0(a;;y)={}_1F_0(a;;x+y-yx).
\eqno\eq{34}
$$
We now give two $q$-analogues of \eqtag{34}, valid in the algebra $\Cqixy$.

\Prop{7}
In the algebra $\Cqixy$ (so $xy=qyx$) we have for $a\in\CC$ the identities
$$
\eqalignno{
{}_1\phi_0(a;;q,x)\,{}_1\phi_0(a;;q,y)&=
{}_1\phi_0(a;;q,x+y-yx),&\eq{35}
\cr
{}_1\phi_0(a;;q,y)\,{}_1\phi_0(a;;q,x)&=
{}_1\phi_0(a;;q,x+y-ayx).&\eq{36}
\cr}
$$

\Proof
In view of \eqtag{37} the identities \eqtag{35}, \eqtag{36} can
be equivalently written as
$$
\eqalignno{
e_q(x)\,E_q(-ax)\,E_q(-ay)\,e_q(y)&=
E_q(-a(x+y-yx))\,e_q(x+y-yx),
\cr
E_q(-ay)\,e_q(y)\,e_q(x)\,E_q(-ax)&=
e_q(x+y-ayx)\,E_q(-a(x+y-ayx)).
\cr}
$$
In view of \eqtag{38}, \eqtag{18}, \eqtag{12} and \eqtag{39}
these identities are in their turn equivalent to
$$
\eqalignno{
e_q(x)\,E_q(-a(x+y))\,e_q(y)&=E_q(-a(x+y-yx))\,e_q(x)\,e_q(y),
\cr
E_q(-ay)\,e_q(x+y)\,E_q(-ax)&=e_q(x+y-ayx)\,E_q(-ay)\,E_q(-ax).
\cr}
$$
Once more, these identities can be rewritten into equivalent forms:
$$
\eqalignno{
E_q(-a(x+y))&=(e_q(x))^{-1}\,E_q(-a(x+y-yx))\,e_q(x),
\cr
e_q(x+y)&=e_q(ay)\,e_q(x+y-ayx)\,(e_q(ay))^{-1}.
\cr}
$$
By a similar argument as in the proof of Proposition \thtag{3}
these two identities will follow if we can show that
$$
\eqalignno{
x+y&=(e_q(x))^{-1}\,(x+y-yx)\,e_q(x),&\eq{40}
\cr
x+y&=e_q(ay)\,(x+y-ayx)\,(e_q(ay))^{-1}.&\eq{41}
\cr}
$$
As for \eqtag{40}, its \RHS\ can be rewritten as
$$
x+(e_q(x))^{-1}\,y(1-x)\,e_q(x)=
x+e_q(x)^{-1}\,y\,e_q(qx)=
x+e_q(x)^{-1}\,e_q(x)\,y=x+y,
$$
where we used \eqtag{13}.
Formula \eqtag{41} can be proved by a similar argument.\hhalmos

\bPP
Formula \eqtag{36} was given by Faddeev \& Volkov
[11, p.315, multiplication rule for $s(\la,w)$],
see also Kirillov [20, Section 2.5, Exercise 3].
The terminating cases of \eqtag{35}, \eqtag{36} (i.e.\ $a=q^{-n}$,
see Kirillov [20, Section 2.5, Lemma 10]) are:

\Cor{10}
In the algebra $\Cqxy$ we have
for $n\in\Zplus$ the identities
$$
(x;q)_n\,(y;q)_n=(x+y-q^nyx;q)_n,\quad
(y;q)_n\,(x;q)_n=(x+y-yx;q)_n.
\eqno\eq{93}
$$

\Sec{4} {Identities for $q$-exponentials with q-Heisenberg relations}
Proposition \thtag{3} can be generalized to an identity in the
algebra $\Cqhxy$ of formal power series
$$
\sum_{k,l,m=0}^\iy a_{k,l,m}\,c^my^lx^k
$$
with arbitrary complex coefficients $a_{k,l,m}$ and where $x,y,c$
satisfy the {\sl q-Heisenberg relations}
$$
xy-qyx=(1-q)c,\quad xc=cx,\quad yc=cy.
\eqno\eq{21}
$$
Here $q\in(0,1)$ is fixed as before.
Note that $c$ is a central element of $\Cqhxy$. Note that,
by adding the relation $c=0$, we can map the algebra $\Cqhxy$ homomorphically
onto the algebra $\Cqixy$.
It is also interesting to remark that the algebra $\Cqhxy$ is isomorphic
to the algebra of formal power series in $x,y,z$ with relations
$$
xy-yx=(1-q)z,\quad
xz-qzx=0,\quad
zy-qyz=0.
\eqno\eq{95}
$$
Just let $z$ and $c$ be related by
$$
z=c-yx.
\eqno\eq{97}
$$

We need the following identity in $\Cqhxy$, which is proved by
induction \wrt $n$:
$$
x^ny=q^nyx^n+(1-q^n)cx^{n-1}.
\eqno\eq{26}
$$

The generalization below of Proposition \thtag{3}
was first given by
Rogers [30] (in operational form, see Bowman [6] for a modern
treatment). The same result was found later, independently, by
Gelfand \& Fairlie [17, (46)],
McDermott \& Solomon [28, (10)] and
Kashaev [18, (4.19)].

\Prop{4}
In the algebra $\Cqhxy$ (i.e. with relations \eqtag{21}) we have the identities
$$
\eqalignno{
e_q(x)\,e_q(y)
&=e_q(y-yx+c)\,e_q(x)&\eq{23}
\cr
&=e_q(y)\,e_q(-yx+c)\,e_q(x)&\eq{24}
\cr
&=e_q(y)\,e_q(x-yx+c).&\eq{25}
\cr}
$$

\Proof
As in the previous proof we have equality \eqtag{17}. By \eqtag{26} we see that
$$
e_q(x)\,y
=\sum_{n=0}^\iy {x^ny\over(q;q)_n}
=\sum_{n=0}^\iy {q^nyx^n\over (q;q)_n}+
\sum_{n=1}^\iy{(1-q^n)cx^{n-1}\over(q;q)_n}
=y\,e_q(qx)+c\,e_q(x).
$$
Hence, by \eqtag{13},
$$
e_q(x)\,y=(y-yx+c)\,e_q(x).
\eqno\eq{27}
$$
Now \eqtag{23} follows from \eqtag{17} and \eqtag{27}.
The other two identities follow by application of
Proposition \thtag{2}.\hhalmos

\Remark{18}
Because of \eqtag{95} and \eqtag{97}, formula \eqtag{24} can be rewritten as
$$
e_q(x)\,e_q(y)=e_q(y)\,e_q\bigl((1-q)^{-1}[x,y]\bigr)\,e_q(x).
\eqno\eq{28}
$$
After rescaling we have, still as an identity in $\Cqhxy$:
$$
e_q\bigl((1-q)x\bigr)\,e_q\bigl((1-q)y\bigr)=
e_q\bigl((1-q)y\bigr)\,e_q\bigl((1-q)[x,y]\bigr)\,e_q\bigl((1-q)x\bigr).
\eqno\eq{100}
$$
Replace $x,y,c$ in \eqtag{21} by $X,Y,(1-q)^{-1}C$, respectively, and
let next $q\uparrow 1$ in \eqtag{21} and \eqtag{100}.
Then we obtain the identity
$$
e^X e^Y=e^Ye^Ce^X\quad\hbox{or equivalently}\quad
e^X e^Y=e^{Y+C}e^X
\eqno\eq{91}
$$
in the algebra of formal power series in $X,Y,C$ satisfying the
{\sl Heisenberg relations}
$$
[X,Y]=C,\quad
[X,C]=0,\quad
[Y,C]=0.
\eqno\eq{96}
$$
Note that the second identity in \eqtag{91} immediately follows from
$$
e^X Y e^{-X}=\exp({\rm ad}\,X)\,Y=Y+[X,Y]=Y+C,
$$
where $({\rm ad}\,U)\,V:=[U,V]=UV-VU$.
A slightly deeper identity in the Heisenberg algebra is obtained by applying the
{\sl Baker-Campbell-Hausdorff formula} (see for instance Varadarajan
[35, Section 2.15])
to the case of relations \eqtag{96}. This yields
$$
e^Y e^X=e^{X+Y-\half C}\quad\hbox{or equivalently}\quad
e^{X+Y}=e^Ye^{\half C}e^X.
\eqno\eq{92}
$$
Formula \eqtag{92} has often been observed in literature.
A $q$-analogue of \eqtag{92} in a $q$-Heisenberg algebra (but not precisely
the algebra $\Cqhxy$ with relations \eqtag{21}) was found by Gelfand \&
Fairlie [17, (49)].
In fact, the result follows easily from \eqtag{12}:

\Prop{19}
In the algebra of formal power series in $x,w,z$ under relations
$$
xw-qwx=(1-q)z^2,\quad
xz=qzx,\quad
zw=qwz
\eqno\eq{98}
$$
we have the identity
$$
e_q(x+w)=e_q(w)\,e_{q^2}(z^2)\,e_q(x).
\eqno\eq{99}
$$

\Proof
It follows from relations \eqtag{98} that
$(x-z)(w+z)=q(w+z)(x-z)$. Hence, by repeated application of \eqtag{12}:
$$
e_q(x+w)=e_q(z+w)\,e_q(x-z)=
e_q(w)\,e_q(z)\,e_q(-z)\,e_q(x).
$$
Now the result follows by application of \eqtag{67}.\hhalmos

\Remark{22}
The algebra of formal power series in $x,w,v$ under relations
$$
xw-qwx=(1-q)v,\quad
xv=q^2vx,\quad
vw=q^2wv
\eqno\eq{103}
$$
is isomorphically embedded into the algebra of formal power series in $x,w,z$
under relations \eqtag{98} by adding the relations $v=z^2$ to relations
\eqtag{103}. This is seen by observing that the algebra of polynomials in
$x,w,v$ under relations \eqtag{103} has a basis of elements $w^kx^lv^m$
($k,l,m\in\Zplus$) and that the algebra of polynomials in
$x,w,z$ under relations \eqtag{98} has a basis of elements $w^kx^lz^m$
($k,l,m\in\Zplus$) (use Bergman's [4] diamond lemma).
Thus it follows from Proposition \thtag{19} that, in the algebra of formal power
series in $x,w,v$ under relations \eqtag{103}, we have the identity
$$
e_q(x+w)=e_q(w)\,e_{q^2}(v)\,e_q(x).
\eqno\eq{104}
$$

\Remark{20}
In relations \eqtag{103} replace $x,w,v$ by
$(1-q)X,(1-q)Y,(1-q)C$ and let
$q\uparrow 1$. Then $X,Y,C$ will satisfy the Heisenberg relations \eqtag{96}
and identity \eqtag{104} becomes the second identity in \eqtag{92}.

We may also rewrite \eqtag{99} as an identity in $\Cqhxy$.
Just put $w:=y[x,y]$ in \eqtag{95}. Then $x,w,z$ satisfy relations \eqtag{98}.
Thus \eqtag{99} becomes in terms of $x,y,z$:
$$
e_q(x+y[x,y])=e_q(y[x,y])\,e_{q^2}\bigl((1-q)^{-1}[x,y]^2\bigr)\,e_q(x).
\eqno\eq{102}
$$
After rescaling we have, still under relations \eqtag{95}:
$$
e_q\bigl((1-q)(x+y[x,y])\bigr)=e_q\bigl((1-q)y[x,y]\bigr)\,
e_{q^2}\bigl((1-q)[x,y]^2\bigr)\,e_q\bigl((1-q)x\bigr).
\eqno\eq{101}
$$
Now make substitutions of $x,y,z$ into $X,Y,C$ as we earlier did for
\eqtag{100}. Next let
$q\uparrow 1$ in \eqtag{21} and \eqtag{101}. Then we obtain the identity
$$
e^{X+YC}=e^{YC}\,e^{\half C^2}\,e^X
$$
under Heisenberg relations \eqtag{96}. This identity
is equivalent to \eqtag{92}.

\Remark{21}
It seems somewhat arbitrary that we stipulate relations \eqtag{98}
in order to obtain identity \eqtag{99}.
In fact, they arise from a much more general Ansatz.
First rewrite \eqtag{104} equivalently as
$$
E_q(-w)\,e_q(x+w)\,E_q(-x)=e_{q^2}(v).
$$
Now we ask more generally under which minimal set of relations
for $x$ and $w$ we have that
$$
E_q(-w)\,e_q(x+w)\,E_q(-x)=f(v)
\eqno\eq{105}
$$
for some formal power series $f$ in one variable and some element
$v$ which is homogeneous of degree 2 in $x$ and $w$, i.e., a linear
combination of $x^2,w^2,xw,wx$.
Surprisingly, the answer is that $v=(1-q)^{-1}(xw-qwx)$,
$f=e_{q^2}$ and that the relations are
$$
x(xw-qwx)=q^2(xw-qwx)x,\quad
(xw-qwx)y=q^2y(xw-qwx),
\eqno\eq{106}
$$
so we recover \eqtag{104} under relations
\eqtag{103} as the unique solution of our problem.
In fact, expansion of the \LHS\ of \eqtag{105} up to quadratic terms
yields $1+(xw-qwx)/(q;q)_2$. It follows from our Ansatz that all
terms of homogeneous odd degree in the expansion of the \LHS\ of \eqtag{105}
must vanish. The vanishing of the third degree terms precisely yields
the relations \eqtag{106}.

\Sec{5} {Equivalence between identities in the non-commuting and the
commuting case}
In the previous sections we saw many examples of identities in
non-commuting variables. It is sometimes
possible to rewrite these identities
in terms of commuting variables, usually in various different ways.
The idea is always the following.
Suppose our identity in non-commuting variables lives in a certain
algebra $\FSA$.
Let $\pi$ be a representation of the algebra
$\FSA$ on a suitable linear
space $\FSF$
of functions (for instance the space of polynomials or formal
power series in one complex variable).
Suppose that there is a subset $\{f_m\}$ of $\FSF$ such that, for all
$a\in\FSA$, we have the implication:
$\pi(a)\,f_m=0\;\hbox{for all $m$}\;\Imply\; a=0$.
An identity $a=b$ in $\FSA$ is then equivalent to the collection of
identities 
$\pi(a)\,f_m=\pi(b)\,f_m$ for all $m$.

As an example, fix $q\in(0,1)$, let $x,y$ be subject to the relation
$xy=qyx$, and
let $\FSA$ be the algebra $\Cqxy$ of polynomials in $x,y$
(cf.\ Section 2)
or the algebra $\Cqixy$ of formal power series in $x,y$
(cf.\ Section 3).
Let $\Cqxy$ resp.\ $\Cqixy$ act on the
space $\FSF$ of polynomials resp.\ formal power series in one
complex variable $z$ by an algebra representation $\pi$
such that
$$
(\pi(x)f)(z):=qz\,f(qz),\quad
(\pi(y)f)(z):=z\,f(z).
\eqno\eq{47}
$$
Then $\pi(x)\,\pi(y)=q\,\pi(y)\,\pi(x)$, so $\pi$ preserves
the relation $xy=qyx$.

By induction with respect to $k$ we see that
$$
(\pi(x^k)f)(z)=q^{\half\,k(k+1)}\,z^k\,f(q^kz)\quad(k\in\Zplus),
$$
so
$$
(\pi(y^lx^k)f)(z)=z^l\,(\pi(x^k)f)(z)=
q^{\half\,k(k+1)}\,z^{k+l}\,f(q^kz)\quad(k,l\in\Zplus).
$$
Thus, by a little abuse of notation,
$$
\pi(y^lx^k)\,z^m=q^{\half k(k+1)+km}\,z^{k+l+m}\quad(k,l,m\in\Zplus).
\eqno\eq{43}
$$
Now we let a formal power series in $x,y$ act on a formal power series
in $z$ in a $\si$-additive way:
$$
\eqalignno{
\pi\left(\sum_{k,l}c_{k,l}\,y^lx^k\right)\,\sum_{m}a_mz^m&=
\sum_{k,l,m}c_{k,l}\,a_m\,q^{\half k(k+1)+km}\,z^{k+l+m}
\cr
&=\sum_{n=0}^\iy z^n
\sum_{k+l+m=n}c_{k,l}\,a_m\,q^{\half k(k+1)+km}.
\cr}
$$
We see that the result is again a well-defined formal power series in $z$.

The reader is warned that a $\si$-additive extension of a \rep\ as we
gave above, is not always possible. Then one has to work with \rep s
by unbounded operators on a Hilbert space, and identities satisfied
by formal series may no longer hold in the \rep, see for instance
Woronowicz [36].
(I thank S. L. Woronowicz for pointing this out to me.)

Next we show that, if $a=\sum_{k,l}c_{k,l}\,y^lx^k$ and
$\pi(a)\,z^m=0$ for all $m\in\Zplus$, then $a=0$.
Indeed, if
$$
\sum_{k,l}c_{k,l}\,\pi(y^lx^k)\,z^m=0\quad\hbox{for all $m\in\Zplus$}
$$
then
$$
\sum_{k,l}c_{k,l}\,q^{\half k(k+1)+km}\,z^{k+l}=0\quad
\hbox{for all $m\in\Zplus$.}
$$
Hence for all $n\in\Zplus$ we have
$$
\sum_{k=0}^n c_{k,n-k}\,q^{\half k(k+1)}\,(q^m)^k=0\quad
\hbox{for all $m\in\Zplus$.}
$$
This shows that the $c_{k,n-k}$ are 0.

Now we will see how we can rewrite the identities \eqtag{3} and \eqtag{12}
involving noncommuting variables into equivalent commutative form by means
of the representation $\pi$ of $\Cqixy$ (or $\Cqxy$).
For fixed $n\in\Zplus$
the $q$-binomial formula \eqtag{3} is equivalent to the set of identities
$$
\pi\bigl((x+y)^n\bigr)\,z^m=\sum_{k=0}^n \qchoose nkq
\pi(y^{n-k} x^k)\,z^m\quad(m\in\Zplus).
\eqno\eq{42}
$$
By induction with respect to $n$ we see that
$$
\pi\bigl((x+y)^n\bigr)\,z^m=(-q^{m+1};q)_n\,z^{m+n}.
\eqno\eq{59}
$$
By also using \eqtag{43} and \eqtag{45} we see that \eqtag{42}
is equivalent to
$$
(-q^{m+1};q)_n\,z^{m+n}=
\sum_{k=0}^n{(q^{-n};q)_k\over(q;q)_k} (-q^{m+n+1})^k\,z^{m+n}\quad
(m\in\Zplus).
$$
We can divide both sides of this last identity by $z^{m+n}$.
Thus we have shown that \eqtag{3} is equivalent to the
terminating $q$-binomial sum \eqtag{9}
considered for all $z=-q^{n+m+1}$ ($m\in\Zplus$),
which in its turn is equivalent to \eqtag{9} considered
for all $z\in\CC$.

Let us next handle identity \eqtag{12} in this way.
Observe first that $(\pi(y)f)(z)=zf(z)$ implies that
$$
\bigl(\pi(g(y))\,f\bigr)(z)=g(z)\,f(z)
\eqno\eq{46}
$$
for any two formal power series $f$ and $g$ in one variable.
Identity \eqtag{12} is equivalent to the set of identities
$$
\pi\left(e_q(x+y)\right)\,z^m=\pi\bigl(e_q(y)\,e_q(x)\bigr)\,z^m
\quad(m\in\Zplus).
$$
Now use \eqtag{46} and expand $e_q(x+y)$ and $e_q(x)$.
We obtain that \eqtag{12} is equivalent to
$$
\sum_{n=0}^\iy{\pi\left((x+y)^n\right)\,z^m\over(q;q)_n}=
e_q(z)\sum_{k=0}^\iy{\pi(x^k) z^m\over(q;q)_k}\quad(m\in\Zplus).
$$
By \eqtag{59} and \eqtag{43} this can be rewritten as
$$
\sum_{n=0}^\iy {(-q^{m+1};q)_n\,z^{m+n}\over(q;q)_n}=
e_q(z)\sum_{k=0}^\iy{q^{\half k(k+1)+km}\,z^{k+m}\over(q;q)_k}\quad
(m\in\Zplus).
$$
On dividing both sides by $z^m$ and on using the first equality in
\eqtag{11} we see that \eqtag{12} is equivalent to
the series of identities
$$
{}_1\phi_0(-q^{m+1};;q,z)=e_q(z)\,E_q(q^{m+1}z)\quad(m\in\Zplus),
$$
i.e., to special cases of the evaluation of non-terminating
$q$-hypergeometric series, cf.\ \eqtag{37}.

I invite the reader to experiment with several other representations $\pi$
of algebras $\Cqxy$ or $\Cqixy$ on the space of polynomials or formal power
series in $z$, for instance
$$
\eqalignno{
&\bigl(\pi(x)f\bigr)(z):=f(qz),\quad
\bigl(\pi(y)f\bigr)(z):=z\,f(z),&\eq{48}
\cr
&\bigl(\pi(x)f\bigr)(z):=(D_qf)(z),\quad
\bigl(\pi(y)f\bigr)(z):=f(qz),&\eq{49}
\cr
&\bigl(\pi(x)f\bigr)(z):=f(z+\log q),\quad
\bigl(\pi(y)f\bigr)(z):=e^z\,f(z),&\eq{50}
\cr
&\bigl(\pi(x)f\bigr)(z):=e^z\,f(z),\quad
\bigl(\pi(y)f\bigr)(z):=f(z-\log q).&\eq{114}
\cr}
$$
In \eqtag{49} we used the notation for the {\sl q-derivative}
$$
(D_qf)(z):={f(z)-f(qz)\over(1-q)z}\,.
\eqno\eq{110}
$$
For instance, consider the $q$-binomial formula \eqtag{3} in representation
\eqtag{48}, with both sides acting on $e_q(z)$, in order to arrive at a
special case of the $q$-Chu-Vandermonde sum [16, (II.6)]
(one upper parameter zero).
Also consider identity \eqtag{19} in representation \eqtag{50}
or identity \eqtag{115} in representation \eqtag{114}, with both
sides acting on functions $e^{i\mu z}$, in order to arrive at an identity
in commuting variables which is equivalent to the summation formula
\eqtag{33} of non-terminating $q$-binomial series. Here the infinite sum
results from the action of the left-hand sides of \eqtag{19} or \eqtag{115}
on $e^{i\mu z}$.
Faddeev \& Kashaev [10, Section 2] point out an alternative
for the action of these left-hand sides. They observe, in representation
\eqtag{114}, that the action of $g_1(y)g_2(x)$ ($g_1$ and $g_2$ suitable
functions)
on $e^{i\mu z}$ can be written as multiplication by a certain double integral
involving $g_1$ and $g_2$ (the symbol of the product of the operators
$\pi(g_1(y))$ and $\pi(g_2(y))$). This argument looks quite formal and the
convergence of the resulting double integral is not clear.

A wealth of further results in the spirit of this section is contained in
Cigler [7], [8].
He also gives applications to continuous $q$-Hermite polynomials and to
$q$-Laguerre polynomials and he develops a $q$-analogue of Rota's
umbral calculus [31].

A well-known \rep\ $\pi$ of the algebra $\Cqhxy$ (see \eqtag{21}, \eqtag{95},
\eqtag{97}) on the space of formal power series is given by
$\pi(x):=(1-q)\,c\,D_q$, $(\pi(y)f)(t):=t\,f(t)$.
See already Rogers [30] for many results using this representation,
or a modern treatment by Bowman [6].

\Sec{6} {The q-logarithm}
{\sl Euler's dilogarithm} is defined by
$$
Li_2(z):=\sum_{n=1}^\iy {z^n\over n^2}\,,\quad |z|<1.
$$
A. N. Kirillov [20, (2.52)]
defines the following $q$-analogue ($0<q<1$ as before):
$$
Li_2(z;q):=\sum_{n=1}^\iy {z^n\over n(1-q^n)}\,,\quad|z|<1.
\eqno\eq{29}
$$
Formally we have the termwise limit
$$
\lim_{q\uparrow 1}(1-q)\,Li_2(z;q)=Li_2(z).
\eqno\eq{116}
$$

Kirillov [20, Section 2.5, Lemma 8] observes the following remarkable formula:
$$
Li_2(z;q)=\log\bigl(e_q(z)\bigr),\quad|z|<1.
\eqno\eq{89}
$$
For the proof note that
$$
Li_2(qz;q)=Li_2(z;q)+\log(1-z)
$$
(substitute the corresponding power series).
Hence
$$
\exp\bigl(Li_2(z;q)\bigr)=
{\exp\bigl(Li_2(qz;q)\bigr)\over 1-z}=\ldots=
{\exp\bigl(Li_2(q^kz;q)\bigr)\over (z;q)_k}.
$$
On taking limits for $k\to\iy$, the \RHS\ tends to
$$
{\exp\bigl(Li_2(0;q)\bigr)\over (z;q)_\iy}=\exp(0)\,e_q(z)=e_q(z).
$$

A more precise formulation of \eqtag{116} going back to
Ramanujan (see Berndt [5, Ch.~27, Entry 6])
gives an
asymptotic series for
$Li_2(z;q)$ (or $\log(e_q(z))$
in rising powers of $-\log q$ as $q\uparrow 1$.
Kirillov [20, Section 2.5, Corollary 10] and Ueno \& Nishizawa [34]
derive this asymptotic series
by using Euler-Maclaurin's summation formula.

Faddeev \& Kashaev [10, S2], see also Kirillov
[20, Theorem I], indicate that Rogers' five-term
identity for Euler's dilogarithm can be obtained as a limit case as $q\uparrow1$
of \eqtag{19} or \eqtag{115}. This uses \eqtag{89} and \eqtag{116}.
Faddeev and Kashaev also use the representation \eqtag{114}
of the \LHS\ of \eqtag{115} by means of a double integral
(see the end of Section 5). Their arguments are quite formal.

Next we consider a $q$-analogue of
$$
-\log(1-z)=\sum_{n=1}^\iy {z^n\over n}\,,\quad|z|<1.
$$
We define and notate it as
$$
\log_q(z):=\sum_{n=1}^\iy{z^n\over 1-q^n}
\eqno\eq{30}
$$
and we consider it either as a convergent power series
for $|z|<1$ or as a formal power series.
Note that we have formally the termwise limit
$$
\lim_{q\uparrow 1}\,(1-q)\,\log_q(z)=-\log(1-z).
\eqno\eq{107}
$$

It follows from \eqtag{29} and \eqtag{30} that
$$
\log_q(z)=z\,Li_2'(z;q).
$$
Hence, by \eqtag{89}:
$$
\log_q(z)={z\,e_q'(z)\over e_q(z)}\,.
\eqno\eq{108}
$$
Another interesting formula is
$$
\log_q(z)=-{\pa\over\pa a}\,{}_1\phi_0(a;;q,z)\Bigr|_{a=1}.
\eqno\eq{109}
$$
It is the $q$-analogue of
$$
-\log(1-z)={\pa\over\pa a} (1-z)^{-a}\Bigr|_{a=0}.
$$
For the proof of \eqtag{109} note that
$$
{\pa\over\pa a}\left({(az;q)_\iy\over(z;q)_\iy}\right)=
{\pa\over\pa a}\left({e_q(z)\over e_q(az)}\right)=
-\,{z\,e_q(z)\,e_q'(az)\over e_q(az)^2}\,.
$$
Hence
$$
{\pa\over\pa a}\left({(az;q)_\iy\over(z;q)_\iy}\right)\Biggr|_{a=1}=
-\,{z\,e_q'(z)\over e_q(z)}\,.
$$
Now apply \eqtag{33} and \eqtag{108}.

Note also the $q$-derivative (see \eqtag{110}) of $\log_q$:
$$
(1-q)\,(D_q \log_q)(z)=\sum_{n=0}^\iy z^n={1\over 1-z}\,.
\eqno\eq{111}
$$
The formulas \eqtag{89}, \eqtag{108}, \eqtag{109}
and \eqtag{111} are of hybrid nature
because they use classical objects (the logarithm, the classical derivative
and the function $z\mapsto (1-z)^{-1}$, respectively) in a context of
$q$-special functions.

It is a natural question to look for the inverse function to $\log_q$,
analogous to the function $x\mapsto 1-e^{-x}$ being inverse to the
function $y\mapsto -\log(1-y)$.
However, the inverse function to $\log_q$ does not seem to have a nice
explicit expression.
(I thank C. Krattenthaler for checking this by use of Maple.)$\;$
Some alternative way to find a $q$-analogue is as follows.
By the chain rule $f(g(x))=x$ implies $f'(g(x))g'(x)=1$.
Now the following $q$-analogue holds:
$$
(D_qf)(g(x))\,(D_qg)(x)=1,\quad
g(x):=1-e_q(-(1-q)x),\quad
f(y):=(1-q)\,\log_q(y).
\eqno\eq{112}
$$
Note however that $(D_qg)(f(y))\,(D_qf)(y)\ne1$ for $f$ and $g$ as in
\eqtag{112}.

The following Proposition gives a $q$-analogue of the classical
functional equation
$$
\log(xy)=\log x+\log y,
\hbox{ or equivalently }
-\log\bigl((1-x)(1-y)\bigr)=-\log(1-x)-\log(1-y),
$$
for commuting variables $x,y$. This Proposition was independently
found by A. N. Kirillov [20, Section 2.5, Exercise 11] and by the author.

\Prop{5}
In the algebra $\Cqixy$ we have the identity
$$
\log_q(x+y-yx)=\log_q(x)+\log_q(y).
\eqno\eq{31}
$$

\Proof
Since $(x-yx)y=qy(x-yx)$ we get by the $q$-binomial formula
\eqtag{3} that
$$
\eqalignno{
(x+y-yx)^n&=\sum_{k=0}^n\qchoose nkq y^{n-k}(x-yx)^k
\cr
&=\sum_{k=0}^n y^{n-k}\,(y;q)_k\,x^k.
\cr}
$$
Hence
$$
\eqalignno{
\log_q(x+y-yx)&=
\sum_{n=1}^\iy\left(
{1\over 1-q^n}\,y^n+
\sum_{k=1}^n{1\over1-q^n}\,
\qchoose nkq\,y^{n-k}\,(y;q)_k\,x^k\right)
\cr
&=\log_q(y)+
\sum_{k=1}^\iy
\left(\sum_{n=k}^\iy{1\over 1-q^n}\,\qchoose nkq\,y^{n-k}\right)
(y;q)_k\,x^k.
\cr}
$$
So we are done if we can show that
$$
\left(\sum_{n=k}^\iy{1-q^k\over 1-q^n}\,\qchoose nkq\,y^{n-k}\right)
(y;q)_k=1.
\eqno\eq{32}
$$
It is sufficient to prove \eqtag{32} for complex $y$ with $|y|<1$.
Then
$$
\eqalignno{
&\sum_{n=k}^\iy{1-q^k\over 1-q^n}\,\qchoose nkq\,y^{n-k}=
\sum_{m=0}^\iy{1-q^k\over 1-q^{m+k}}\,\qchoose{m+k}kq\,y^m
\cr
&\quad
=\sum_{m=0}^\iy{1-q^k\over 1-q^{m+k}}\,{(q^{k+1};q)_m\over(q;q)_m}\,y^m
=\sum_{m=0}^\iy{(q^k;q)_m\over(q;q)_m}\,y^m=
{(q^ky;q)_\iy\over(y;q)_\iy}={1\over(y;q)_k}\,,
\cr}
$$
where we used in the forelast identity the evaluation \eqtag{33}
of a non-terminating
$q$-binomial series.\hhalmos

\bPP
A somewhat formal, but very quick proof of Proposition \thtag{5} is obtained
from \eqtag{35}. Just differentiate both sides of \eqtag{35} with respect to
$a$, put $a=1$ and use \eqtag{109}.

Yet another proof is obtained from \eqtag{12}, which we use in the form
$e_q(t(x+y))=e_q(ty)\,e_q(tx)$ with $t\in\RR$.
Differentiate both sides with respct to $t$, put $t=1$ and use \eqtag{108}
in order to obtain
$$
e_q(x+y)\,\log_q(x+y)=e_q(y)\,(\log_q(y)+\log_q(x))\,e_q(x).
$$
Replace next $e_q(x+y)$ by $e_q(y)\,e_q(x)$ and pull the $e_q(x)$ factor
through $\log_q(x+y)$, on using \eqtag{113}.

\Sec{7} {Jackson integral in q-commuting variables}
This section reviews results from Kempf \& Majid [19, Section 1].
We start with a $q$-analogue of Taylor series for
$q$-commuting variables.

Fix $q\in(0,1)$. The {\sl q-derivative} $D_qf$ of a function $f$ on $\RR$
was already defined in \eqtag{110}.
If $f(x)=x^n$ then $(D_qf)(x)=\displaystyle{1-q^n\over 1-q}\,x^{n-1}$,
so we can let $D_q$ and its iterates $D_q^k$ act on polynomials or
formal power series in $x$. We have
$$
D_q^k x^n={(q^n;q^{-1})_k\over(1-q)^k}\,x^{n-k}.
$$
Hence, in the algebra $\Cqxy$
we can rewrite the $q$-binomial formula \eqtag{3} as
$$
(x+y)^n=\sum_{k=0}^n{(1-q)^k\over(q;q)_k}\,y^k\,D_q^k x^n.
\eqno\eq{51}
$$
Now let $f(x):=\sum_{n=0}^\iy c_nx^n$ be a formal power series.
Then, in the algebra $\Cqixy$ of formal power series in $x,y$ under relation
$xy=qyx$, we have
$$
\eqalignno{
&f(x+y)=
\sum_{n=0}^\iy c_n(x+y)^n
=\sum_{n=0}^\iy\sum_{k=0}^n c_n{(1-q)^k\over(q;q)_k}\,y^k\,D_q^k x^n
\cr
&=\sum_{k=0}^\iy
{(1-q)^k\over(q;q)_k}\,y^k\sum_{n=k}^\iy c_n D_q^kx^n
=\sum_{k=0}^\iy {(1-q)^k\over(q;q)_k}\,y^k\,D_q^k\,f(x).
\cr}
$$
Thus we have proved:

\Prop{11}
Let $f$ be a formal power series in one variable.
Let $xy=qyx$. Then, in the algebra $\Cqixy$ we have for each $m\in\Zplus$
the identity
$$
\eqalignno{
f(x+y)&=
\sum_{k=0}^\iy {1\over(q;q)_k}\,y^k\,((1-q)D_q)^k\,f(x)&\eq{117}
\cr
&=\sum_{k=0}^{m-1}{1\over(q;q)_k}\,y^k\,((1-q)D_q)^k\,f(x)+
y^m\,g_m(x,y)&\eq{52}
\cr}
$$
for a suitable element $g_m(x,y)$ of $\Cqixy$.

\bPP
If we write $\FSO(y^m)$ instead of $y^m\,g_m(x,y)$ then \eqtag{52}
implies in particular, for $m=2$, that
$$
f(x+y)=f(x)+y\,D_q\,f(x)+\FSO(y^2)\quad(xy=qyx),
\eqno\eq{53}
$$
as a $q$-analogue in $q$-commuting variables of the classical formula
$$
f(x+y)=f(x)+yf'(x)+\FSO(y^2)\quad(xy=yx).
$$

Fix $q\in(0,1)$.
For a function $f$ on $\RR$, the {\sl Jackson integral} is defined by
$$
\int_0^x f(t)\,d_qt:=(1-q)\sum_{k=0}^\iy f(q^kx)\,q^kx,
\eqno\eq{56}
$$
where $x\in\RR$, provided the sum on the \RHS\ converges absolutely, for
instance if $f$ is bounded near zero.
The Jackson integral avant la lettre of $f(t):=t^n$ was already
computed by Fermat:
$$
\int_0^x t^n\,d_qt={1-q\over1-q^{n+1}}\,x^{n+1}.
\eqno\eq{54}
$$
Thus, if
$$
f(x):=\sum_{n=0}^\iy c_n x^n
\eqno\eq{55}
$$
is a formal power series in
a formal (not necessarily real or complex) variable $x$ then
we obtain its Jackson integral
$$
\int_0^x f(t)\,d_qt=\sum_{n=0}^\iy c_n\,{1-q\over1-q^{n+1}}\,x^{n+1}
\eqno\eq{60}
$$
also as a formal power series.

Let now $xy=qyx$. Then we obtain in the algebra $\Cqxy$ by \eqtag{54}
and a twofold appication of the $q$-binomial formula \eqtag{3}:
$$
\eqalignno{
\int_0^{x+y} t^n\,d_qt
&={1-q\over1-q^{n+1}}\,(x+y)^{n+1}
\cr
&={1-q\over1-q^{n+1}}\,y^{n+1}+
{1-q\over1-q^{n+1}}\sum_{l=0}^n\qchoose {n+1}lq y^l x^{n+1-l}
\cr
&={1-q\over1-q^{n+1}}\,y^{n+1}+(1-q)
\sum_{l=0}^n\qchoose nlq {1\over 1-q^{n-l+1}}\,y^lx^{n-l+1}
\cr
&={1-q\over1-q^{n+1}}\,y^{n+1}+(1-q)
\sum_{l=0}^n\sum_{k=0}^\iy\qchoose nlq y^l(q^kx)^{n-l+1}
\cr
&={1-q\over1-q^{n+1}}\,y^{n+1}+(1-q)
\sum_{k=0}^\iy(q^kx+y)^n\,q^kx
\cr
&=\int_0^y t^n\,d_qt+\int_0^x(t+y)^n\,d_qt,
\cr}
$$
where we use the definition
$$
\int_0^x f(t+y)\,d_qt:=(1-q)\sum_{k=0}^\iy f(q^kx+y)\,q^kx.
$$
Thus we have shown that
$$
\int_0^{x+y} t^n\,d_qt=
\int_0^y t^n\,d_qt+\int_0^x(t+y)^n\,d_qt\quad(xy=qyx)
$$
and therefore also for formal power series \eqtag{55}:
$$
\int_0^{x+y} f(t)\,d_qt=
\int_0^y f(t)\,d_qt+\int_0^x f(t+y)\,d_qt\quad(xy=qyx).
\eqno\eq{58}
$$
Now recall the definition
$$
\int_a^b f(t)\,d_qt:=\int_0^a f(t)\,d_qt-\int_0^b f(t)\,d_qt,
\eqno\eq{57}
$$
where the two Jackson integrals on the \RHS\ are defined by \eqtag{56}.
Definition \eqtag{57} can also be used for formal variables $a,b$.
Thus by \eqtag{58} and \eqtag{57} we have the following Proposition:

\Prop{12}
Let $f$ be a formal power series in one variable.
Let $xy=qyx$. Then, in the algebra $\Cqixy$ we have
$$
\int_y^{x+y} f(t)\,d_qt=\int_0^x f(t+y)\,d_qt.
\eqno\eq{62}
$$

\bPP
Thus the translation invariance of the Riemann integral, which
seemed to be destroyed when $q$-deforming it to the Jackson integral,
can be preserved when we work with $q$-commuting variables.

By way of example we give a fourth proof of the functional equation
\eqtag{31} for the $q$-logarithm in $q$-commuting variables.
First observe, by straightforward application of \eqtag{60} and
\eqtag{30}, that
$$
\log_q(x)=\int_0^x (1-t)^{-1}\,d_qt.
\eqno\eq{61}
$$
Now let us work in $\Cqixy$. Then
$$
\eqalignno{
&\log_q(x+y-xy;q)-\log(y;q)
=\int_0^{x+y-yx} (1-t)^{-1}\,d_qt-
\int_0^y(1-t)^{-1}\,d_qt
\cr
&=\int_0^{x-yx}(1-(t+y))^{-1}\,d_qt
=(1-q)\,\sum_{k=0}^\iy
\bigl(1-q^k(x-yx)-y\bigr)^{-1}\,q^k(x-yx)
\cr
&=(1-q)\,\sum_{k=0}^\iy
\bigl((1-y)(1-q^kx)\bigr)^{-1}\,(1-y)q^kx
=(1-q)\,\sum_{k=0}^\iy(1-q^kx)^{-1}\,q^kx\
\cr
&=\int_0^x(1-t)^{-1}\,d_qt
=\log_q(x),
\cr}
$$
where we applied \eqtag{62} in the second equality.
Thus we have given a new proof of \eqtag{31}.

\Sec{8} {q-Hermite polynomials and a q-Fourier transform pair} 
The {\sl discrete q-Hermite I polynomials}
(see Koekoek \& Swarttouw [22, Section 3.28] and references given there)
are given by
$$
\eqalignno{
h_n(x;q):=&x^n\,{}_2\phi_0(q^{-n},q^{-n+1};;q^2,q^{2n-1}x^{-2})
\cr
=&(q;q)_n\,\sum_{k=0}^{[n/2]}{(-1)^k\,q^{k(k-1)}\,x^{n-2k}\over
(q^2;q^2)_k\,(q;q)_{n-2k}}\,.&\eq{147}
\cr}
$$
They are orthogonal polynomials satisfying the orthogonality relations
$$
\int_{-1}^1 h_m(x;q)\,h_n(x;q)\,E_{q^2}(-q^2x^2)\,d_qx=
b_q\,q^{\half n(n-1)}\,(q;q)_n\,\de_{m,n},
\eqno\eq{139}
$$
where
$$
b_q:=(1-q)\,(q,-q,-1;q)_\iy.
\eqno\eq{158}
$$
There is the following generating function:
$$
{(t^2;q^2)_\iy\over (xt;q)_\iy}=E_{q^2}(-t^2)\,e_q(xt)=
\sum_{n=0}^\iy {h_n(x;q)\over (q;q)_n}\,t^n
\quad(|xt|<1).
\eqno\eq{138}
$$
Several useful formulas can be derived from this generating function.
First we can expand a monomial:
$$
x^n=(q;q)_n\,\sum_{k=0}^{[\half n]}{h_{n-2k}(x;q)\over
(q^2;q^2)_k\,(q;q)_{n-2k}}.
\eqno\eq{137}
$$
For the proof, multiply both sides of \eqtag{138} with $e_{q^2}(t^2)$
and next compare coefficients of $t$.

Next we have
$$
\sum_{k=0}^m{(q^{-m};q)_k\over(q;q)_k}\,q^k\,h_k(x;q)\,x^{m-k}=
\cases{
(-1)^n\,q^{-n^2}(q;q^2)_n&if $m=2n$,
\cr
0& if $m=2n+1$.
\cr}
\eqno\eq{77}
$$
For the proof, multiply both sides of \eqtag{138} with $E_q(-xt)$
and next compare coefficients of $t$.

>From the generating function \eqtag{138} together with the orthogonality
relations \eqtag{139} we obtain
$$
\int_{-1}^1 e_q(-ixt)\,h_n(x;q)\,E_{q^2}(-q^2x^2)\,d_qx=
b_q\,q^{\half n(n-1)}\,i^{-n}\,t^n\,e_{q^2}(-t^2).
\eqno\eq{140}
$$
For the proof, replace $t$ by $-it$ in \eqtag{138}, multiply both sides
with $h_n(x;q)\,E_{q^2}(-q^2x^2)\,e_{q^2}(-t^2)$ and $q$-integrate both sides
over $[-1,1]$.

The {\sl discrete q-Hermite II polynomials}
(see Koekoek \& Swarttouw [22, Section 3.29] and references given there)
are given by
$$
\eqalignno{
\wt h_n(x;q):=&x^n\,{}_2\phi_1(q^{-n},q^{-n+1};0;q^2,-q^2 x^{-2})&\eq{150}
\cr
=&(q;q)_n\,\sum_{k=0}^{[\half n]}
{(-1)^k\,q^{-2nk}\,q^{k(2k+1)}\,x^{n-2k}\over(q^2;q^2)_k\,(q;q)_{n-2k}}\,.
&\eq{143}
\cr}
$$
(There is a slight error in [22, (3.29.1)], which is corrected
in formula \eqtag{150} above.)
They are related to the discrete $q$-Hermite I polynomials by
$$
h_n(ix;q^{-1})=i^n\,\wt h_n(x;q).
\eqno\eq{142}
$$
They are orthogonal polynomials satisfying orthogonality relations
given by a Jackson integral over $(-\iy,\iy)$.
Let us use the notation
$$
\int_0^{\ga.\iy}f(t)\,d_qt:=(1-q)\sum_{k=-\iy}^\iy f(q^k\ga)\,q^k\ga
\eqno\eq{63}
$$
for the Jackson integral over $(0,\iy)$ of a function defined on
$\{q^k\ga\mid k\in\ZZ\}$
for some $\ga\in(0,\iy)$, where we suppose that the sum on the right
hand side of \eqtag{63} absolutely converges. So the definition is dependent
on $\ga$ but invariant under the transform $\ga\mapsto q\ga$.
For $f$ defined on
$\{\pm q^k\ga\mid k\in\ZZ\}$
we can also define the Jackson integral of $f$
over $(-\iy,\iy)$, again depending on $\ga$ and invariant under
$\ga\mapsto q\ga$:
$$
\eqalignno{
\int_{-\ga.\iy}^{\ga.\iy}f(t)\,d_qt:=&
\left(\int_0^{\ga.\iy}f(t)\,d_qt-\int_0^{-\ga.\iy}f(t)\,d_qt\right)&
\eq{64}
\cr
=&(1-q)\sum_{k=-\iy}^\iy\bigl(f(q^k\ga)+f(-q^k\ga)\bigr)q^k\ga.&\eq{118}
\cr}
$$
We again suppose that the infinite sums are absolutely convergent.

The orthogonality relations for the discrete $q$-Hermite II polynomials
are:
$$
\int_{-\ga.\iy}^{\ga.\iy}
\wt h_m(x;q)\,\wt h_n(x;q)\,e_{q^2}(-x^2)\,d_qx=
c_q(\ga)\,q^{-n^2}\,(q;q)_n\,\de_{m,n},
\eqno\eq{151}
$$
where
$$
c_q(\ga):=
{2(1-q)\,(q^2,-q\ga^2,-q\ga^{-2};q^2)_\iy\,\ga\over
(-\ga^2,-q^2/\ga^2,q;q^2)_\iy}\,.
\eqno\eq{157}
$$
Note that the orthogonality measure is not uniquely determined.
There is the generating function
$$
{(-xt;q)_\iy\over(-t^2;q^2)_\iy}=
e_{q^2}(-t^2)\,E_q(xt)=
\sum_{n=0}^\iy {q^{\half n(n-1)}\over(q;q)_n}\,\wt h_n(x;q)\,t^n.
\eqno\eq{141}
$$
Now we can obtain formulas analogous to \eqtag{137}--\eqtag{140},
either by using \eqtag{141} just as \eqtag{138} was used in the discrete
$q$-Hermite I case, or by using \eqtag{142}:
$$
x^n=(q;q)_n\sum_{k=0}^{[\half n]}{q^{-2nk+3k^2}\over(q;q)_{n-2k}\,(q^2;q^2)_k}\,
\wt h_{n-2k}(x;q),
\eqno\eq{144}
$$
$$
\sum_{k=0}^m{(q^{-m};q)_k\over(q;q)_k}\,q^{mk}\,\wt h_k(x;q)\,x^{m-k}=
\cases{(-1)^n\,(q;q^2)_n&if $m=2n$,
\cr
0&if $m=2n+1$,
\cr}
\eqno\eq{145}
$$
$$
\int_{-\ga.\iy}^{\ga.\iy}
E_q(iqxt)\,\wt h_n(x;q)\,e_{q^2}(-x^2)\,d_qx=
c_q(\ga)\,q^{-\half n(n-1)}\,i^n\,t^n\,E_{q^2}(-q^2t^2).
\eqno\eq{146}
$$

Now we can combine formulas for cases I and II of the discrete $q$-Hermite
polynomials.
It follows from \eqtag{137}, \eqtag{140} and \eqtag{143} that
$$
\int_{-1}^1 e_q(-ixt)\,x^n\,E_{q^2}(-q^2x^2)\,d_qx=
b_q\,q^{\half n(n-1)}\,i^{-n}\,\wt h_n(t;q)\,e_{q^2}(-t^2)
\eqno\eq{148}
$$
and it follows from \eqtag{144}, \eqtag{146} and \eqtag{147} that
$$
\int_{-\ga.\iy}^{\ga.\iy}
E_q(iqxt)\,x^n\,e_{q^2}(-x^2)\,d_qx=
c_q(\ga)\,q^{-\half n(n-1)}\,i^n\,h_n(t;q)\,E_{q^2}(-q^2t^2).
\eqno\eq{149}
$$
By comparing \eqtag{140} with \eqtag{149} or \eqtag{148} with \eqtag{146},
we arrive at the following pair of $q$-Fourier transforms, which are inverse
to each other when acting on suitable functions:

\Theor{26}
Let
$$
\eqalignno{
(\FSF_qf)(y)&:={1\over b_q}\int_{-1}^1 e_q(-ixy)\,f(x)\,d_qx,&\eq{171}
\cr
(\wt\FSF_{q,\ga}g)(x)&:={1\over c_q(\ga)}
\int_{-\ga.\iy}^{\ga.\iy}E_q(iqxy)\,g(y)\,d_qy.&\eq{172}
\cr}
$$
Then the following two statements are equivalent:
\item{(a)}
$f(x)={\rm polynomial}(x)\times E_{q^2}(-q^2x^2)$ and $g=\FSF_q f$;
\item{(b)}
$g(y)={\rm polynomial}(y)\times e_{q^2}(-y^2)$ and $f=\wt\FSF_{q,\ga} g$.
\sLP
In particular:
$$
\eqalignno{
f(x)=h_n(x;q)\,E_{q^2}(-q^2x^2)\;&\Longleftrightarrow\;
g(y)=q^{\half n(n-1)}\,i^{-n}\,y^n\,e_{q^2}(-y^2);&\eq{153}
\cr
f(x)=x^n\,E_{q^2}(-q^2x^2)\;&\Longleftrightarrow\;
g(y)=q^{\half n(n-1)}\,i^{-n}\,\wt h_n(y;q)\,e_{q^2}(-y^2).&\eq{154}
\cr}
$$

\bPP
By $q$-integration by parts and by using \eqtag{13} we can see how $\FSF_q$ and
$\wt\FSF_{q,\ga}$ send a $q$-derivative operator to a multiplication operator
and a multiplication operator to a $q$-derivative operator.
For this purpose we also need a variant $D_q^+$,
called {\sl forward q-derivative},
of the ({\sl backward})
$q$-derivative $D_q=D_q^-$ as defined in \eqtag{110}:
$$
(D_q^-f)(x):={f(x)-f(qx)\over (1-q)x},\quad
(D_q^+f)(x):={f(q^{-1}x)-f(x)\over (1-q)x}.
\eqno\eq{152}
$$

\Prop{27}
If $f$ is continuous at 0 and $f(q^{-1})=0=f(-q^{-1})$ then
$$
(1-q)\,(\FSF_q(D_q^+ f))(y)=iy\,(\FSF_q f)(y).
$$
If $g$ is continuous at 0 and
$\lim_{n\to\iy}E_q(\pm ix\ga q^{-n+1})\,g(\pm\ga q^{-n})=0$ then
$$
(1-q)\,(\wt\FSF_{q,\ga}(D_q^- g))(x)=-ix\,(\wt\FSF_{q,\ga}g)(x).
$$

\bPP
It follows from \eqtag{153}, \eqtag{154} and Proposition \thtag{27}
that (with $D_q^\pm$ acting on functions of $x$):
$$
\eqalignno{
(1-q)\,D_q^+ \bigl(h_n(x;q)\,E_{q^2}(-q^2x^2)\bigr)&=
-q^{-n}\,h_{n+1}(x;q)\,E_{q^2}(-q^2x^2),&\eq{155}
\cr
(1-q)\,D_q^- \bigl(\,\wt h_n(x;q)\,e_{q^2}(-x^2)\bigr)&=
-q^n\,\wt h_{n+1}(x;q)\,e_{q^2}(-x^2).&\eq{156}
\cr}
$$
These two formulas can also be proved independently. 
Thus, if \eqtag{155} and \eqtag{156} are given then the general case of
\eqtag{153} and \eqtag{154} follows from the special case $n=0$ by means
of Proposition \thtag{27}.

Note that iteration of \eqtag{155}, \eqtag{156} yields Rodrigues type
formulas
$$
\eqalignno{
h_n(x;q)&=(-1)^n\,q^{\half n(n-1)}\,e_{q^2}(q^2x^2)\,
(1-q)^n\,(D_q^+)^n E_{q^2}(-q^2x^2),&\eq{160}
\cr
\wt h_n(x;q)&=(-1)^n\,q^{-\half n(n-1)}\,E_{q^2}(x^2)\,
(1-q)^n\,(D_q^-)^n\,e_{q^2}(-x^2).&\eq{161}
\cr}
$$

All formulas, given in this section, involving discrete $q$-Hermite
polynomials
until now are analogues of formulas for classical Hermite polynomials:
just take limits as $q\uparrow 1$ (after possibly some rescaling).
However, there are some other nice formulas for classical Hermite polynomials
for which the $q$-analogue can be better given with $q$-commuting variables.
For instance, in the algebra $\Cqxy$ ($xy=qyx$) we have:
$$
h_n(x+y;q)=\sum_{k=0}^n \qchoose nkq y^{n-k}\,h_k(x;q)
\eqno\eq{163}
$$
For the proof, let $t$ be scalar and use the generating function \eqtag{138}
and the functional equation \eqtag{12}:
$$
\eqalignno{
&\sum_{n=0}^\iy {h_n(x+y;q)\over (q;q)_n}\,t^n=
e_q((x+y)t)\,E_{q^2}(-t^2)=
e_q(yt)\,e_q(xt)\,E_{q^2}(-t^2)
\cr
&=\sum_{l=0}^\iy{y^l\over(q;q)_l}\,t^l\,
\sum_{k=0}^\iy{h_k(x;q)\over (q;q)_k}\,t^k
=\sum_{n=0}^\iy\left(\sum_{k=0}^n \qchoose nkq y^{n-k}\,h_k(x;q)\right)
{t^n\over(q;q)_n}.
\cr}
$$

Another example of a formula for Hermite polynomials without $q$-analogue
in commuting variables is the formula expanding $H_n(\la x)$ in terms of
a series of $H_m(x)$:
$$
H_n(\la x)=n!\,\sum_{k=0}^{[\half n]}
{(-1)^k\,(1-\la^2)^k\,\la^{n-2k}\over (n-2k)!\,k!}\,H_{n-2k}(x).
\eqno\eq{178}
$$
This formula can be proved by use of the generating function
$$
e^{2xt-t^2}=\sum_{n=0}^\iy {H_n(x)\,t^n\over n!}\,.
$$
Just write
$$
e^{2\la xt-t^2}=e^{2\la xt-\la^2t^2}\,e^{-(1-\la^2)t^2},
$$
expand all factors in terms of powers of $t$ and compare coefficients of
$t^n$ on both
sides.
Note that we used that
$e^{-t^2}=e^{-\la^2 t^2}\,e^{-(1-\la^2)t^2}$, which strongly suggests to use
$q$-commuting variables for a $q$-analogue.

Suppose now that $\la,\mu$ satisfy the relation $\la\mu=q^\half\mu\la$ and let
$x$ and $t$ be scalar.
Then $\la^2\mu^2=q^2\mu^2\la^2$, so, by \eqtag{38}:
$$
e_q(\la xt)\,E_{q^2}(-(\la^2+\mu^2)t^2)=
\bigl(e_q(\la xt)\,E_{q^2}(-\la^2t^2)\bigr) E_{q^2}(-\mu^2t^2).
$$
Now expand factors in terms of powers of $t$ by using the generating function
\eqtag{138} and compare coefficients of $t^n$ on both sides.
Then we obtain a $q$-analogue of \eqtag{178}:
$$
\sum_{k=0}^{[\half n]}{(-1)^k\,q^{k(k-1)}\,\la^{n-2k}\,(\la^2+\mu^2)^k\,
x^{n-2k}\over(q;q)_{n-2k}\,(q^2;q^2)_k}=
\sum_{k=0}^{[\half n]}
{(-1)^k\,q^{k(k-1)}\,\la^{n-2k}\,\mu^{2k}\,h_{n-2k}(x;q)\over
(q;q)_{n-2k}\,(q^2;q^2)_k}\,.
\eqno\eq{179}
$$
The \LHS\ of \eqtag{179}, after multiplication by $(q;q)_n$,
can be considered as $h_n(x;q)$ being kind of rescaled by means of $\la$ and
$\mu$.
Formula \eqtag{179} can be used in order to arrive at a $q$-analogue
of the Fourier transform sending
$H_n(x)\,e^{-\half x^2}$ to a constant multiple of this, or more generally
$H_n(x)\,e^{-ax^2}$ to a constant multiple of
$H_n((4a(1-a))^{-\half} x)\,e^{-x^2/(4a)}$ (see [9, (1.10.8),
(2.10.10)]).
Just combine \eqtag{149} with \eqtag{179}. I do not give details, but the
reader should compare with very related results in
Finkelstein \& Marcus [12], where the quantum group $SU_q(2)$
is also brought into the game.

The $q$-Fourier transform $\wt \FSF_{q,\ga}$ (see \eqtag{172}) also occurs in
Kempf \& Majid [19, Section VIB]. However, the Fourier kernel is
written there as $e_q(iyx)$ with $x$ and $y$ $q$-commuting. 
Their discussion is tied up very much with the notion of translation invariant
Jackson integral over $(-\iy,\iy)$ and with $\CC_q[x]$ considered as a
braided Hopf algebra, see the next two sections for a few more details.

As far as I know, a purely analytic approach to the $q$-Fourier transform pair
in Theorem \thtag{26} has not earlier appeared in literature. Of course,
there are many natural further questions, e.g.\ extension of the transforms to
a bigger class of functions and continuity properties of the transforms.

\Sec{9} {Translation invariance of Jackson integral over $(-\iy,\iy)$}
The results of this section are essentially due to Kempf \& Majid
[19, Section IV], but the approach is different.

In \eqtag{64}, \eqtag{118} we gave the definition of a Jackson integral
over $(-\iy,\iy)$. We cannot extend this definition to the case
$\int_{-x.\iy}^{x.\iy}f(t)\,d_qt$ where $f$ is a formal power series and
$x$ is a formal variable (compare with \eqtag{60} for the Jackson integral from
0 to $x$), since the Jackson integral of $f(t):=t^n$ over
$(-\iy,\iy)$ is not well-defined.

Still, in view of the classical formula
$$
\int_{-\iy}^\iy f(t)\,dt=\int_{-\iy}^\iy f(t+y)\,dt\quad(y\in\RR),
$$
valid for absolutely convergent integrals, we would like to find an
extension of the Jackson integral translation invariance \eqtag{62}
on finite intervals for $q$-commuting variables
to the case of an infinite interval.
So we would like to have that
$$
\int_{-x.\iy}^{x.\iy} f(t)\,d_qt=
\int_{-x.\iy}^{x.\iy} f(t+y)\,d_qt\quad(xy=qyx)
\eqno\eq{65}
$$
for suitable formal power series $f$ for which both
sides of \eqtag{65} have meaning.

Let me first give a completely formal proof of \eqtag{65}, as a limit case
of \eqtag{62}. For $r\in\Zplus$ it follows from \eqtag{62} that
$$
\int_{-q^{-r}(x+q^ry)}^{q^{-r}(x+q^ry)}f(t)\,d_qt=
\int_{-q^{-r}x}^{q^{-r}x} f(t+y)\,d_qt\quad(xy=qyx).
$$
Now let $r\to\iy$. Then equality \eqtag{65} is obtained as a formal limit
case.

Evidently this argument gives only heuristic evidence for the validity of
\eqtag{65}.
Let me give next a still formal, but more satisfactory proof of
\eqtag{65} for suitable functions $f$.
By formal substitution of \eqtag{118} in the \RHS\ of \eqtag{65} and by
\eqtag{117} we have
$$
\eqalignno{
\int_{-x.\iy}^{x.\iy} f(t+y)\,d_qt
&=(1-q)\sum_{k=-\iy}^\iy\bigl(f(q^kx+y)+f(-q^kx+y)\bigr)q^kx
\cr
&=\sum_{m=0}^\iy y^m\,(1-q)\sum_{k=-\iy}^\iy
\bigl(f_m(q^kx)+f_m(-q^kx)\bigr)\,q^kx
\cr
&=\sum_{m=0}^\iy y^m\int_{-x.\iy}^{x.\iy} f_m(t)\,d_qt,&\eq{119}
\cr}
$$
where
$$
f_m(z)={1\over(q;q)_m}\,((1-q)D_q)^m\,f(z).
\eqno\eq{122}
$$
If $f$ is a function on $\RR$ then so is $f_m$.
So the Jackson integrals $\int_{-x.\iy}^{x.\iy}f_m(t)\,d_qt$
will have concrete meaning for $x\in\RR$ and if the sums defining the
Jackson integral are convergent.
However, $x$ must $q$-commute with $y$, so $x$ cannot be in $\RR$.
We can circumvent this dilemma by passing to a suitable representation
of the relation $xy=qyx$.
Let us take a slight extension of the representation \eqtag{48},
now using the dot notation instead of $\pi$:
$$
x.g(z):=\ga\,g(qz),\quad
y.g(z):=z\,g(z).
\eqno\eq{120}
$$
Here $\ga\in\RR\backslash\{0\}$ is fixed and $g(z)$ is a formal power series.
Thus
$$
x.z^k=\ga q^k z^k\quad(k\in\Zplus),
$$
which we will formally extend to
$$
f(x).z^k=f(\ga q^k)\,z^k\quad(k\in\Zplus)
$$
if $f$ is a function on $\RR$.
Hence
$$
\left(\int_{-x.\iy}^{x.\iy} f(t)\,d_qt\right).\,z^k=
\left(\int_{-\ga q^k.\iy}^{\ga q^k.\iy} f(t)\,d_qt\right)z^k=
\left(\int_{-\ga.\iy}^{\ga.\iy} f(t)\,d_qt\right)z^k
$$
provided that the sums defining the Jackson integral
$$
I_f(\ga):=\int_{-\ga.\iy}^{\ga.\iy} f(t)\,d_qt
$$
converge absolutely.
Note that $I_f(\ga)=I_f(q\ga)$. For such $f$ we conclude that, for any formal
power series $g$:
$$
\left(\int_{-x.\iy}^{x.\iy} f(t)\,d_qt\right).\,g(z)=I_f(\ga)\,g(z).
\eqno\eq{124}
$$

Now take up \eqtag{119} again in the representation \eqtag{120}.
We find
$$
\left(\int_{-x.\iy}^{x.\iy} f(t+y)\,d_qt\right).\,g(z)=
\sum_{m=0}^\iy I_{f_m}(\ga)\,y^m.g(z)=
\sum_{m=0}^\iy I_{f_m}(\ga)\,z^m\,g(z),
$$
while
$$
\left(\int_{-x.\iy}^{x.\iy} f(t)\,d_qt\right).\,g(z)=I_f(\ga)\,g(z).
$$
So formula \eqtag{65} in the \rep\ \eqtag{120}
is equivalent with the vanishing of $I_{f_m}(\ga)$ ($m=1,2,\ldots$).
It is easy to find a class of functions $f$ for which these numbers vanish.
Note that each $f_m$ ($m\in\Zplus$) is a $q$-derivative of another function.

\Lemma{23}
Let $\ga\in\RR\backslash\{0\}$ and let $f$ be a function on
$\{\pm\ga q^k\mid k\in\ZZ\}$ such that
$\lim_{k\to\iy} f(q^k\ga)=\lim_{k\to\iy}f(-q^k\ga)$ and
$\lim_{k\to\iy} f(\pm q^{-k}\ga)=0$.
Then $\int_{-\ga.\iy}^{\ga.\iy}(D_qf)(t)\,d_qt=0$.

\Proof
It follows by summation by parts that
$$
\int_{-\ga.\iy}^{\ga.\iy}(D_qf)(t)\,d_qt
=\lim_{m,n\to\iy}\bigl(f(q^{-m}\ga)-f(q^n\ga)-f(-q^{-m}\ga)+f(-q^n \ga)\bigr).
\eqno\halmos
$$

\Prop{24}
Let $\ga\in\RR\backslash\{0\}$ and let $f$ be a function on
$\{\pm\ga q^k\mid k\in\ZZ\}\cup \{0\}$ such that, for all $m\in\Zplus$,
$D_q^mf$ is continuous at 0 and
$$
|(D_q^m f)(\pm q^{-k}\ga)|=\FSO(q^{(1+\ep)k})\quad \hbox{as $k\to\iy$}
\eqno\eq{123}
$$
for certain $\ep>0$. Then $f$ satisfies the translation invariance \eqtag{65}
in the \rep\ \eqtag{120}.
If, moreover, the estimate \eqtag{123} is satisfied for all $\ep>0$
(so $f$ and all its $q$-derivatives are {\sl rapidly decreasing} on the domain
of definition),
then $f$ multiplied with any polynomial also 
satisfies the translation invariance \eqtag{65}
in the \rep\ \eqtag{120}.

\Proof
Because of \eqtag{122} and the estimate on $D_q^mf$, the Jackson integrals
defining the $I_{f_m}(\ga)$ converge absolutely.
For each $m>0$, $f_m$ is a $q$-derivative,
so $I_{f_m}(\ga)=0$ by Lemma \thtag{23} and the estimate for $D_q^{m-1}f$.
For the proof of the last statement use the $q$-Leibniz rule
$(D_q(fg))(x)=f(x)\,(D_qg)(x)+(D_qf)(x)\,g(qx)$.\hhalmos

\bPP
Note that the continuity at 0 in the Proposition is satisfied if $f$ is
the restriction of a function which is $C^\iy$ on a neighbourhood of 0.

As an example we consider Jackson integrals involving the $q$-Gaussian
$g_q(x)=e_{q^2}(-x^2)$ (cf.\ \eqtag{66}).
This function satisfies the conditions of the
Proposition (including the rapid decreasing property).
In fact, it follows by induction \wrt $m$ that
$(D_q^m g_q)(x)=p_m(x)\,e_{q^2}(-x^2)$ with $p_m$ a polynomial of degree $m$
(more concretely a discrete $q$-Hermite II polynomial, see
\eqtag{161}).
Now observe that
$$
|e_{q^2}(-x^2)|=\prod_{k=0}^\iy|1+q^{2k}x^2|^{-1}\le(1+q^{2n}x^2)^{-n}=
\FSO(|x|^{-2n})\quad\hbox{as $x\to\pm\iy$}
$$
for all $n\in\Zplus$.
Thus we know that \eqtag{65} in the \rep\ \eqtag{120}
is valid for $f(x):=x^m\,e_{q^2}(-x^2)$.
Let us see the implications of this result for $q$-special functions,
by which we can make contact with the results of Section 8.

It follows from \eqtag{147} that
$$
h_m(0;q)=\cases{
(-1)^n\,q^{n(n-1)}\,(q;q^2)_n&if $m=2n$,
\cr
0&if $m=2n+1$.
\cr}
\eqno\eq{174}
$$
Hence the case $t=0$ of \eqtag{149} yields
$$
\int_{-\ga.\iy}^{\ga.\iy} t^m\,e_{q^2}(-t^2)\,d_qt=\cases{
c_q(\ga)\,q^{-n^2}(q;q^2)_n&if $m=2n$,
\cr
0&if $m=2n+1$
\cr}
\eqno\eq{69}
$$
where $c_q(\ga)$ is given by \eqtag{157}.
In fact, the case $m=2n$ of \eqtag{69} is a consequence of Ramanujan's
${}_1\psi_1$ summation formula, see [16, (II.29)].

Consider \eqtag{65} with $f(x):=x^m\,e_{q^2}(-x^2)$.
The \LHS\ of \eqtag{65}, when acting as an operator on a formal
power series $g(z)$ in the \rep\
\eqtag{120},
can be evaluated as
$$
\left(\int_{-x.\iy}^{x.\iy} t^m\,e_{q^2}(-t^2)\,d_qt\right).\,g(z)=\cases{
c_q(\ga)\,q^{-n^2}(q;q^2)_n\,g(z)&if $m=2n$,
\cr
0&if $m=2n+1$
\cr}
\eqno\eq{75}
$$
We expand the \RHS\ of \eqtag{65}, still formally, as
$$
\eqalignno{
&\int_{-x.\iy}^{x.\iy} e_{q^2}(-(t+y)^2)\,(t+y)^m\,d_qt
\cr&\quad
=\sum_{k=0}^m \qchoose mkq e_{q^2}(-y^2)
\int_{-x.\iy}^{x.\iy} e_q(-yt)\,e_{q^2}(-t^2)\,y^{m-k}\,t^k\,d_qt
\cr&\quad
=\sum_{k=0}^m q^{-(k+1)(m-k)} \qchoose mkq e_{q^2}(-y^2)
\left(\int_{-x.\iy}^{x.\iy} e_q(-yt)\,e_{q^2}(-t^2)\,t^k\,d_qt\right)
y^{m-k},&\eq{72}
\cr}
$$
where we used \eqtag{3} and \eqtag{68} for the first equality,
while the second
equality follows from \eqtag{118} and the $q$-commutation of
$x$ and $y$. Now we give meaning to
the \RHS\ of \eqtag{72} as an
operator acting on a formal power series $g(z)$ in the \rep\ \eqtag{120}.
Consider first:
$$
\eqalignno{
&\left(\int_{-x.\iy}^{x.\iy} e_q(-yt)\,e_{q^2}(-t^2)\,t^k\,d_qt\right).\,g(z)
\cr
&=\sum_{l=0}^\iy{(-1)^l q^{\half l(l-1)}y^l\over (q;q)_l}\,
\left(\int_{-x.\iy}^{x.\iy} e_{q^2}(-t^2)\,t^{k+l}\,d_qt\right).\,g(z)
\cr
&=\sum_{l=0}^\iy{(-1)^l q^{\half l(l-1)}z^l\over (q;q)_l}\,
\left(\int_{-\ga.\iy}^{\ga.\iy} e_{q^2}(-t^2)\,t^{k+l}\,d_qt\right)g(z)
\cr
&=\left(\int_{-\ga.\iy}^{\ga.\iy} E_q(-zt)\,e_{q^2}(-t^2)\,t^k\,d_qt\right)g(z)
\cr
&=c_q(\ga)\,q^{-\half k(k-1)}\,i^k\,h_k(iq^{-1}z;q)\,E_{q^2}(z^2)\,g(z).
&\eq{}
\cr}
$$
Here we used \eqtag{124} and \eqtag{149}.
After substitution of this result in \eqtag{72} we obtain that the \RHS\
of \eqtag{72} acting on $g(z)$ becomes:
$$
c_q(\ga)\,\sum_{k=0}^m {(q^{-m};q)_k\over(q;q)_k}\,q^k\,i^{-m}\,
h_k(iq^{-1}z;q)\,(iq^{-1}z)^{m-k}\,g(z).
\eqno\eq{159}
$$
We know, at least formally, that
\eqtag{75} and \eqtag{159} must be equal to each other.
But this equality can equivalently be written as
\eqtag{77}, which we had already proved in an elementary way.

Let us next consider whether \eqtag{65} holds when $f$ equals
the second $q$-Gaussian (cf.\ \eqtag{66})
$$
G_q(x):=E_{q^2}(-x^2)=(x^2;q^2)_\iy.
$$
Note that $G_q(\pm q^{-m})=0$ for $m\in\Zplus$.
So $G_q(\pm\ga q^{-m})=0$ for $m$ a sufficiently large integer
if $\ga$ is an integer
power of $q$. However, if $\ga$ is not an integer power of $q$ then
$|G_q(\ga q^{-m})|$ increases faster than $C^m$ for any $C>1$ as $m\to\iy$.
Indeed, take $n\in\ZZ$ such that
$\ga^2 q^{-2n-2}-1\ge C$. Then, for $m\ge n$:
$$
|G_q(\ga q^{-m})|\ge
(\ga^2q^{-2n-2}-1)^{m-n}\,|G(\ga q^{-n})|.
$$
So the Jackson integral $\int_{-\ga.\iy}^{\ga.\iy} G_q(t)\,t^m\,d_qt$
(the analogue of \eqtag{69}) only converges absolutely
for $\pm\ga$ being an integer
power of $q$ and then it turns down to computing the Jackson integral over
$[-q,q]$.

It follows from \eqtag{143} that
$$
\wt h_m(0;q)=\cases{
(q;q^2)_n\,(-1)^n\,q^{n-2n^2}&if $m=2n$,
\cr
0&if $m=2n+1$.
\cr}
\eqno\eq{169}
$$
Hence the case $t=0$ of \eqtag{148} yields
$$
\int_{-q}^q E_{q^2}(-t^2)\,t^m\,d_qt=
\cases{b_q\,q^{2n+1}\,(q;q^2)_n&if $m=2n$,
\cr
0&if $m=2n+1$,
\cr}
\eqno\eq{126}
$$
where $b_q$ is given by \eqtag{158}.
Alternatively,
formula \eqtag{126} can also be obtained by a completely elementary
computation.

Since $G_q$ is a $C^\iy$-function and since it vanishes on the set
$\{\pm q^{-m}\mid m\in\Zplus\}$, it clearly satisfies all conditions of
Proposition \thtag{24} for $\ga=1$. So the function $f(x):=x^m\,E_{q^2}(-x^2)$
will satisfy \eqtag{65} in the \rep\ \eqtag{120} for $\ga=1$, i.e.\ in the
\rep\ \eqtag{48}. In order to see the implications of this, we can imitate
what we did for the other $q$-Gaussian $g_q$.
The \LHS\ of \eqtag{65}, when acting as an operator on a formal
power series $g(z)$ in the \rep\
\eqtag{48},
can be evaluated as
$$
\left(\int_{-x.\iy}^{x.\iy} t^m\,E_{q^2}(-t^2)\,d_qt\right).\,g(z)=\cases{
b_q\,q^{2n+1}(q;q^2)_n\,g(z)&if $m=2n$,
\cr
0&if $m=2n+1$
\cr}
\eqno\eq{127}
$$
by using \eqtag{126}.
We expand the \RHS\ of \eqtag{65}, still formally, as
$$
\eqalignno{
&\int_{-x.\iy}^{x.\iy} (t+y)^m\,E_{q^2}(-(t+y)^2)\,d_qt
\cr&\quad
=\sum_{k=0}^m \qchoose mkq y^{m-k}
\int_{-x.\iy}^{x.\iy} t^k\,E_{q^2}(-t^2)\,E_q(-yt)\,
E_{q^2}(-y^2)\,d_qt
\cr&\quad
=\sum_{k=0}^m \qchoose mkq y^{m-k}
\left(\int_{-x.\iy}^{x.\iy} t^k\,E_{q^2}(-t^2)\,E_q(-yt)\,
d_qt\right)E_{q^2}(-q^{-2}y^2).&\eq{128}
\cr}
$$
Here we used \eqtag{3} and \eqtag{133}.
Now we give meaning to
the \RHS\ of \eqtag{128} as an
operator acting on a formal power series $g(z)$ in the \rep\ \eqtag{48}.
First derive, analogous to the proof of \eqtag{} but now using \eqtag{148},
that
$$
\eqalignno{
\left(\int_{-x.\iy}^{x.\iy} t^k\,E_{q^2}(-t^2)\,E_q(-yt)\,
d_qt\right).\,g(z)
&=\left(\int_{-q}^q t^k\,E_{q^2}(-t^2)\,e_{q^2}(-q^{-2}zt)\,d_qt\right)g(z)
\cr
&=qb_q\,i^{-k}\,q^{\half k(k+1)}\,e_{q^2}(q^{-2}z^2)\,
\wt h_k(-iq^{-1}z;q)\,g(z).
\cr}
$$
After substitution of this result in \eqtag{128} we obtain that the \RHS\
of \eqtag{128} acting on $g(z)$ becomes
$$
qb_q\sum_{k=0}^m {(q^{-m};q)_k\over(q;q)_k}\,q^{(m+1)k}\,i^k\,z^{m-k}\,
\wt h_k(-iq^{-1}z;q)\,g(z).
\eqno\eq{162}
$$
The \RHS\ of \eqtag{127}
must be equal to \eqtag{162}.
But this equality can equivalently be written as
\eqtag{144}, which we had already proved
in an elementary way.

Thus we have seen in this section that the translation invariance
\eqtag{65} for the case that $f(t)=t^mg_q(t)$ or $t^mG_q(t)$ turns down to
the identities \eqtag{77} and \eqtag{144} for discrete $q$-Hermite polynomials.

\Sec{10} {Braided Hopf algebras}
In this section we introduce braided Hopf algebras and show the relevance of
this structure for the results of Section 8.
First we recall the notion of an ordinary Hopf algebra
(see for instance
Abe [1], Sweedler [33], Koornwinder [25, Section 1]).

We will work over the field of complex numbers, so a {\sl linear space} will
mean a complex linear space.
If $V$ and $W$ are linear spaces then the {\sl tensor product} $V\ten W$
will be the linear space which is the algebraic tensor product of $V$ and $W$,
so $V\ten W$ will be spanned by the elements $v\ten w$ ($v,w\in V$).
The tensor products $V\ten\CC$ and $\CC\ten V$ will be naturally identified
with $V$.

By an {\sl algebra} we will mean a complex
associative algebra with identity element 1.
The field of complex numbers is an algebra in an evident way.
If $\FSA$ is an algebra then the mapping $m\colon\FSA\ten\FSA\to\FSA$
will denote the linear extension of the bilinear mapping
$(a_1,a_2)\mapsto a_1a_2\colon\FSA\times\FSA\to\FSA$.
If $\FSA$ and $\FSB$ are algebras then an {\sl algebra homomorphism}
$\phi\colon\FSA\to\FSB$ is a linear mapping satisfying
$\phi(a_1a_2)=\phi(a_1)\phi(a_2)$ and $\phi(1)=1$.

If $V$ is a linear space then the {\sl flip operator}
is the linear operator $\si\colon V\ten V\to V\ten V$ such that
$\si(v\ten w)=w\ten v$.
If $\FSA$ is an algebra then we give an algebra structure to
$\FSA\ten\FSA$ by putting $(a\ten b)(c\ten d):=ac\ten bd$ and
by extending this to a bilinear mapping of $(\FSA\ten\FSA)\times(\FSA\ten\FSA)$
to $\FSA\ten\FSA$. Let
$m_{\FSA\ten\FSA}\colon \FSA\ten \FSA\ten \FSA\ten \FSA\to\FSA\ten\FSA$
be the linear operator $m$ corresponding to this algebra structure.
Thus $m_{\FSA\ten\FSA}$ acting on $a\ten b\ten c\ten d$
can be written as the composition of two operators:
$$
a\ten b\ten c\ten d\longrightarrow
a\ten c\ten b\ten d\longrightarrow ac\ten bd,
$$
so
$$
m_{\FSA\ten\FSA}=(m_\FSA\ten m_\FSA)\circ(\id\ten\si\ten\id).
\eqno\eq{79}
$$

\Def{17}
A {\sl Hopf algebra} is an algebra $\FSA$ equipped with
three additional operators
$\De\colon\FSA\to\FSA\ten\FSA$ ({\sl comultiplication}),
$\ep\colon\FSA\to\CC$ ({\sl counit}) and
$S\colon\FSA\to\FSA$ ({\sl antipode}), where $\De$ and $\ep$ are algebra
homomorphisms and $S$ is a linear mapping, and where the following additional
properties are satisfied:
$$
\eqalignno{
&(\De\ten\id)\circ\De=(\id\ten\De)\circ\De\quad
({\sl coassociativity}),&\eq{82}
\cr
&(\ep\ten\id)\circ\De=\id=(\id\ten\ep)\circ\De,&\eq{83}
\cr
&\bigl(m\circ(S\ten\id)\circ\De\big)(a)=\ep(a)\,1=
\bigl(m\circ(\id\ten S)\circ\De\big)(a)\quad(a\in\FSA).&\eq{84}
\cr}
$$

\bPP
It can be shown as a consequence of this definition that the antipode
is anti-mul\-ti\-plica\-tive and anti-comultiplicative:
$$
\eqalign{
&S(ab)=S(b)S(a),\quad S(1)=1,
\cr
&(S\ten S)\circ\si\circ\De=\De\circ S,\quad
\ep\circ S=\ep.
\cr}
\eqno\eq{86}
$$

In the definition of Hopf algebra the flip operator $\si$ entered
in the specification \eqtag{79} of the algebra structure of $\FSA\ten\FSA$,
and this algebra structure is needed since $\De\colon\FSA\to\FSA\ten\FSA$
is required to be an algebra homomorphism.
In a {\sl braided Hopf algebra}
the role of $\si$ is taken over by some
other bijective
linear mapping $\Psi\colon\FSA\ten\FSA\to\FSA\ten\FSA$, the so-called
{\sl braiding}. The multiplication in $\FSA\ten\FSA$ is now defined by
$$
(a\ten b)(c\ten d):=(m\ten m)(\id\ten\Psi\ten\id)(a\ten b\ten c\ten d).
\eqno\eq{80}
$$
The braiding $\Psi$ has to satisfy some further axioms such that multiplication
in $\FSA\ten\FSA$ is associative and comultiplication in
$\FSA\ten\FSA$ is coassociative, which I will not give here.
The definition of Hopf algebra is now precisely as in Definition \thtag{17},
but with modified multiplication rule in $\FSA\ten\FSA$.
Braided Hopf algebras were introduced by Majid,
see [27] and references given there.
I refer to his papers for further details.

By way of example consider the {\sl braided line} $\FSA:=\CC_q[x]$
(see Koornwinder [23, Section 6.8], Majid [26]).
As an algebra it is just the algebra $\CC[x]$ of polynomials in one variable
$x$, generated by the element $x$, and with basis $1,x,x^2,\ldots\;$.
Now let $q\in(0,1)$ and introduce the braiding $\Psi$ by specifying it on a
basis of $\FSA\ten\FSA$:
$$
\Psi(x^k\ten x^l):=q^{kl}\,x^l\ten x^k\quad(k,l\in\Zplus).
$$
So for multiplication in $\FSA\ten\FSA$ we will have:
$$
\eqalignno{
(x^k\ten 1)(1\ten x^l)&=x^k\ten x^l,
\cr
(1\ten x^k)(x^l\ten 1)&=q^{kl}\,x^l\ten x^k,
\cr
(x^{k_1}\ten x^{k_2})\,(x^{l_1}\ten x^{l_2})&=
q^{k_2l_1}\,x^{k_1+l_1}\ten x^{k_2+l_2}.
\cr}
$$
Hence
$$
\eqalignno{
\bigl((x^{k_1}\ten x^{k_2})\,(x^{l_1}\ten x^{l_2})\bigr)\,(x^{m_1}\ten x^{m_2})&=
q^{k_2l_1+k_2m_1+l_2m_1}\,x^{k_1+l_1+m_1}\ten x^{k_2+l_2+m_2}
\cr
&=(x^{k_1}\ten x^{k_2})\,\bigl((x^{l_1}\ten x^{l_2})\,
(x^{m_1}\ten x^{m_2})\bigr),
\cr}
$$
which proves the associativity.

Note that $\FSA\ten\FSA$ can be considered as the algebra with generators
$1\ten x$ and $x\ten 1$ and with relation
$(1\ten x)(x\ten 1)=q(x\ten 1)(1\ten x)$, so it is isomorphic with the
algebra $\Cqxy$ under the isomorphisms $x^l\ten x^k\mapsto y^lx^k$.

More generally we can make the $n$-fold tensor product $\ten^n\FSA$
into an algebra by the rule
$$
(x^{k_1}\ten\cdots\ten x^{k_n})(x^{l_1}\ten\cdots\ten x^{l_n})=
q^{\sum_{i>j}k_il_j}\,x^{k_1+l_1}\ten\cdots\ten x^{k_n+l_n}.
$$
A simple computation shows that the multiplication is associative.
Futhermore, the linear subspace of $\ten^n\FSA$ spanned by the elements
$$
1\ten\cdots\ten 1\ten x^{k_1}\ten1\ten\cdots\ten 1\ten x^{k_2}\ten 1\ten\cdots
\cdots\ten 1\ten x^{k_r}\ten 1\ten\cdots\ten 1
$$
(non-zero powers of $x$ only allowed at positions $i_1,i_2,\ldots,i_r$)
is a subalgebra of $\ten^n\FSA$ isomorphic to $\ten^r\FSA$.

Since the comultiplication $\De\colon\FSA\to\FSA\ten\FSA$
has to be an algebra homomorphism, it is sufficient to define it
on the generator $x$ of $\FSA$:
$$
\De(x):=x\ten1+1\ten x.
$$
Then
$$
\bigl((\De\ten\id)\circ\De\bigr)(x)=
x\ten 1\ten 1+1\ten x\ten 1+1\ten 1\ten x=
\bigl((\id\ten\De)\circ\De\bigr)(x),
$$
so the coassociativity, being valid on the generator $x\in\FSA$, will be
valid in general.
If $f$ is a polynomial in one variable then
$$
\De(f(x))=f(x\ten 1+1\ten x).
\eqno\eq{167}
$$
For the comultiplication applied to a general basis element $x^n\in\FSA$
we find
$$
\De(x^n)=\sum_{k=0}^n\qchoose nkq x^{n-k}\ten x^k.
\eqno\eq{81}
$$
Indeed, by rewriting the two sides of \eqtag{81} we have to prove that
$$
(x\ten 1+1\ten x)^n=\sum_{k=0}^n \qchoose nkq (x\ten 1)^{n-k}\,(1\ten x)^k,
$$
and this is true by the $q$-binomial formula \eqtag{3} since
$(1\ten x)(x\ten 1)=q(x\ten 1)(1\ten x)$.

In order to find the counit $\ep\colon\FSA\to\CC$
we appply the first identity of \eqtag{83} to $x^n$ and we use \eqtag{81}:
$$
x^n=(\ep\ten\id)(\De(x^n))=\sum_{k=0}^n\qchoose nkq \ep(x^{n-k})\,x^k,
$$
which yields
$$
\ep(x^n)=\cases{
0&if $n>0$,
\cr
1&if $n=0$.
\cr}
\eqno\eq{87}
$$

Finally we look for the existence of an antipode $S\colon\FSA\to\FSA$.
It turns out that $S(x^n)$ is uniquely found by letting one the identitities
of \eqtag{84} (say the second one) act on the basis elements $x^n$ of
$\FSA$. We obtain
$$
S(x^n):=(-1)^n\,q^{\half n(n-1)}\,x^n.
\eqno\eq{85}
$$
Indeed, the second identity of \eqtag{84} yields
$$
\sum_{k=0}^n\qchoose nkq x^{n-k}\,S(x^k)=
\bigl(m\circ(\id\ten S)\circ\De\bigr)(x^n)=\ep(x^n)\,1=\de_{n,0}\,1,
\eqno\eq{175}
$$
so $S(x^n)$ is obtained by recurrence \wrt $n$.
Now the \LHS\ of \eqtag{175} with
$S(x^k)$ given by \eqtag{85} becomes
$$
\sum_{k=0}^n\qchoose nkq (-1)^k\,q^{\half k(k-1)}\,x^n=
\sum_{k=0}^n{(q^{-n};q)_k\over(q;q)_k}\,q^{nk}\,x^n=
{}_1\phi_0(q^{-n};;q,q^n)\,x^n=(1;q)_n\,x^n,
$$
which equals the \RHS\ of \eqtag{175}.
The proof that the first identity of \eqtag{84} acting on $x^n$ also
holds, amounts to the same computation as we just gave.

Observe that
$$
\eqalignno{
S(x^{m+n})&=q^{mn}\,S(x^m)\,S(x^n),
\cr
\De(S(x^n))&=\sum_{k=0}^n\qchoose nkq q^{k(n-k)}\,S(x^{n-k})\ten S(x^k).
\cr}
$$
Hence
$$
S\circ m=m\circ(S\ten S)\circ\Psi,\quad
\De\circ S=(S\ten S)\circ\Psi\circ\De.
\eqno\eq{176}
$$
This can be considered as an analogue of \eqtag{86} for the braided case,
where the flip $\si$ is replaced by the braiding $\Psi$.

We can reformulate formula \eqtag{163} by using the above
comultiplication:
$$
\De(h_n(x;q))=\sum_{k=0}^n \qchoose nkq x^{n-k}\ten h_k(x;q).
\eqno\eq{168}
$$
When we apply $m\circ(S\ten\id)\circ\De$ to both sides of \eqtag{168} then
we obtain
$$
h_n(0;q)=\sum_{k=0}^n\qchoose nkq (-1)^{n-k}\,q^{\half(n-k)(n-k-1)}\,x^{n-k}\,
h_k(x;q).
$$
In view of \eqtag{174} this is just \eqtag{77}.

We can extend the above braided Hopf algebra structure of $\CC_q[x]$
to the algebra $\CC_q[[x]]$ of formal power series in $x$.
Then $\De$, $\ep$ and $S$ are well-defined in a termwise way on
$\CC_q[[x]]$, by using \eqtag{81}, \eqtag{87} and \eqtag{85}.
For $\De(f(x))$, $f$ being a formal power series, we can still use \eqtag{167}.
In particular, for $\De$, $\ep$ and $S$ acting on $e_q(x)$ we obtain
$$
\De(e_q(x))=e_q(x)\ten e_q(x),\quad
\ep(e_q(x))=1,\quad
S(e_q(x))=E_q(-x).
\eqno\eq{88}
$$
The first identity,
which was already observed in Koornwinder [23, Section 6.8],
follows from (and is equivalent to) \eqtag{12}
since $\De(e_q(x))=e_q(x\ten 1+1\ten x)$ and
$(1\ten x)(x\ten 1)=q(x\ten 1)(1\ten x)$.
Now it follows from \eqtag{84} and \eqtag{88} that
$$
1=\ep(e_q(x))1=\bigl(m\circ(S\ten\id)\circ\De\bigr)(e_q(x))=
S(e_q(x))\,e_q(x)=E_q(-x)\,e_q(x),
$$
i.e., we have reobtained \eqtag{14}.

Formula \eqtag{65} (the translation invariance of the Jackson integral
over $(-x.\iy,x.\iy)$) can also be rephrased in terms of the above
comultiplication. We will work in a very formal way.
Let $\int$ be the linear operator defined by
$$
{\mediumint} (f(x)):=\int_{-x.\iy}^{x.\iy} f(t)\,d_qt.
$$
If $y$ is an element with the property that $xy=q^kyx$ for some integer $k$,
then $y$ commutes with $\int(f(x))$, so $\int (f(x))$ may be considered
as a scalar.
Now we can rewrite \eqtag{65} as
$$
\bigl((\id\ten{\mediumint})\circ\De\bigr)(f(x))=
{\mediumint}(f(x))\;1.
\eqno\eq{170}
$$
Indeed, the \LHS\ can formally be written as
$$
\int_{-x.\iy}^{x.\iy}f(x\ten 1+1\ten t)\,d_qt
$$
and $(1\ten t)(x\ten 1)=q(x\ten 1)(1\ten t)$ if $t=q^kx$ for some integer $k$.

Next we will consider the $q$-Fourier transform $\FSF_q$ defined by \eqtag{171}
from the point of view of this comultiplication.
We fix $q\in(0,1)$ and we rewrite \eqtag{171} very formally as
$$
\FSF_y(f(x)):=\int_{-x.\iy}^{x.\iy} e_q(-ity)\,f(t)\,d_qt\quad(y\in\RR).
$$
Then the linear operator $\FSF_y$, like $\int$, will map to a space of scalars.
It follows from \eqtag{170} and \eqtag{12} that
$$
\eqalignno{
\FSF_y(f(x))\;1
&=(\id\ten\mediumint)\bigl(\De(e_q(-ixy))\,\De(f(x))\bigr)
\cr
&=(\id\ten\mediumint)\Bigl(\bigl(e_q(-ixy)\ten e_q(-ixy)\bigr)\,\De(f(x))\Bigr)
\cr
&=(\id\ten\FSF_y)(\De(f(x))\;e_q(-ixy)
\cr}
$$
Now multiply the \LHS\ and the \RHS\ both by $E_q(ixy)$.
Then we obtain
$$
(\id\ten\FSF_y)(\De(f(x))=\FSF_y(f(x))\;E_q(ixy).
\eqno\eq{173}
$$
This formula is a $q$-analogue of the well-known property of the Fourier
transform that
$$
\int_{-\iy}^\iy e^{-ity}\, f(x+t)\,dt=e^{ixy} \int_{-\iy}^\iy e^{-ity}\,
f(t)\,dt.
$$
Formula \eqtag{173} can also be written as
$$
(\id\ten\mediumint)\bigl((1\ten e_q(-ixy))\;\De(f(x))\bigr)=
(\id\ten\mediumint)\Bigl(\bigl((S\ten\id)\circ\De\bigr)(e_q(-ixy))\,
(1\ten f(x))\Bigr).
$$
More generally we have
$$
(\id\ten\mediumint)\bigl((1\ten g(x))\;\De(f(x))\bigr)
=(\id\ten\mediumint)\Bigl(\bigl((S\ten\id)\circ\De\bigr)(g(x))\;
(1\ten f(x))\Bigr),
\eqno\eq{177}
$$
which is a $q$-analogue of
$$
\int_{-\iy}^\iy g(t)\,f(x+t)\,dt=
\int_{-\iy}^\iy g(t-x)\,f(t)\,dt.
$$
Formula \eqtag{177} is given by Kempf \& Majid in [19, (139)] with
diagrammatic proof in their Figure 2(b).
See this paper also for many further results about Jackson integral
and $q$-Fourier transform in connection with the braided line.

\Sec{11} {Further results}
In this final section I very briefly mention three further results.
\item{1.}
Special function identities involving non-commuting variables
satisfying relations which may be more complicated than $q$-commutation
occur very naturally as addition formulas obtained from quantum groups.
A prototype of this is the addition formula for little
$q$-Legendre polynomials, in the context of $SU_q(2)$.
See Koornwinder [24], where the formula with non-commuting variables
is next equivalently rewritten in commuting variables.
In the context of $U_q(n)$ a similar but more complicated addition formula
in non-commuting variables was obtained by Floris [14] for
$q$-disk polynomials. Next Floris \& Koelink [15] found an
equivalent form in commuting variables of this addition formula.
See Koelink [21] for other examples involving $SU_q(2)$ of
rewriting addition formulas from non-commutative form into commutative form.

\item{2.}
Some $q$-hypergeometric series which are not summable when parameters
and argument
commute, may suddenly become summable when these variables do not commute
but satisfy certain relations. G. Andrews (private communication) showed
me a surprising example of this involving a ${}_{m+1}\phi_m$.
Many further results in this direction, often in operational form,
can be found in Bowman [6].

\item{3.}
The following was communicated to me by A. Yu.\ Volkov.
The second part of formula \eqtag{93} nicely generalizes to
an identity in $x,y,c$ with $q$-Heisenberg relations \eqtag{21}:
$$
(y;q)_n\, (x;q)_n=\prod_{k=0}^{n-1}(1-q^k(x+y-yx+c)+q^{2k}c).
$$
The proof is by induction with respect to $n$. On letting $n\to\iy$ and
taking the inverse on both sides, we obtain an addition to
Proposition \thtag{4}.

\bLP\bLP
{\sl Acknowledgements}\quad
This paper was written while the author was a guest at the
Institute Mittag-Leffler, Djursholm, Sweden during October, November 1995.
I am grateful to the director of the Institute and to the organizers of
the special year ``Analysis on Lie Groups'' for hospitality.
I also thank
G. Andrews,
R. Askey,
D. Bowman,
A. N. Kirillov,
C. Krattenthaler,
A. Yu.\ Volkov
and S. L. Woronowicz
for helpful remarks.

\bLP\bLP
{\bf References}\frenchspacing

\refitem{[1]}
E. Abe,
{\sl Hopf algebras},
Cambridge University Press, 1980.

\refitem{[2]}
R. Askey,
{\sl Abstract of ``Extensions of Hermite polynomials and other orthogonal
polynomials''},
SIAM Newsletter ``Orthogonal Polynomials and Special Functions''
5 (1994), no.1, p.14.

\refitem{[3]}
R. Askey,
{\sl A brief introduction to the world of q},
in {\sl Symmetries and integrability of difference equations},
D. Levi \& P. Winternitz (eds.),
CRM Proceedings and Lecture Notes 9, Amer. Math. Soc., Providence, 1996.

\refitem{[4]}
G. M. Bergman,
{\sl The diamond lemma for ring theory},
Advances in Math. 29 (1978), 178--218.

\refitem{[5]}
B. C. Berndt,
{\sl Ramanujan's notebooks, Part IV},
Springer, 1994.

\refitem{[6]}
D. Bowman,
{\sl q-Difference operators and symmetric expansions},
in preparation.

\refitem{[7]}
J. Cigler,
{\sl Operatormethoden f\"ur $q$-Identit\"aten},
Monatsh. Math. 88 (1979), 87--105.

\refitem{[8]}
J. Cigler,
{\sl Operatormethoden f\"ur $q$-Identit\"aten II: $q$-Laguerre-Polynome},
Monatsh. Math. 91 (1981), 105--117.

\refitem{[9]}
A. Erd\'elyi,
{\sl Tables of integral transforms, Vol. I},
McGraw-Hill, 1954.

\refitem{[10]}
L. Faddeev \& R. M. Kashaev,
{\sl Quantum dilogarithm},
Modern Phys. Lett. A 9 (1994), 427--434.

\refitem{[11]}
L. Faddeev \& A. Yu. Volkov,
{\sl Abelian current algebra and the Virasoro algebra on the lattice},
Phys. Lett. B 315 (1993), 311--318.

\refitem{[12]}
R. Finkelstein \& E. Marcus,
{\sl Transformation theory of the q-oscillator},
J. Math. Phys. 36 (1995), 2652--2672.

\refitem{[13]}
R. Floreanini \& L. Vinet,
{\sl Automorphisms of the q-oscillator algebra and basic orthogonal
polynomials},
Phys. Lett. A 180 (1993), 393--401.

\refitem{[14]}
P. G. A. Floris,
{\sl Addition formula for q-disk polynomials},
Compositio Math. (to appear).

\refitem{[15]}
P. G. A. Floris \& H. T. Koelink,
{\sl A commuting q-analogue of the addition formula for disk polynomials},
Constr. Math. (to appear).

\refitem{[16]}
G. Gasper \& M. Rahman,
{\sl Basic hypergeometric series},
Encyclopedia of Mathematics and its Applications 35,
Cambridge University Press, 1990.

\refitem{[17]}
I. M. Gelfand \& D. B. Fairlie,
{\sl The algebra of Weyl symmetrised polynomials and its quantum extension},
Comm. Math. Phys. 136 (1991), 487--499.

\refitem{[18]}
R. M. Kashaev,
{\sl Heisenberg double and pentagon relation},
preprint, 1995, q-alg/9503005, to appear in Algebraic Analysis.

\refitem{[19]}
A. Kempf \& S. Majid,
{\sl Algebraic $q$-integration and Fourier theory on quantum and braided
spaces},
J. Math. Phys. 35 (1994), 6802--6837.

\refitem{[20]}
A. N. Kirillov,
{\sl Dilogarithm identities},
Progress Theor. Phys. Supplement No. 118 (1995), 61--142.

\refitem{[21]}
H. T. Koelink,
{\sl Addition formulas for q-special functions},
Proceedings ``Special Functions,
q-Series and Related Topics, Toronto, June 1995'',
The Fields Institute Communications Series (to appear).

\refitem{[22]}
R. Koekoek \& R. F. Swarttouw,
{\sl The Askey scheme of hypergeometric orthogonal polynomials and its
$q$-analogue},
Report 94-05, Technical University Delft, Faculty of Technical Mathematics
and Informatics, 1994.

\refitem{[23]}
T. H. Koornwinder,
{\sl Orthogonal polynomials in connection with quantum groups},
in {\sl Orthogonal polynomials: theory and practice},
P. Nevai (ed.),
NATO ASI Series C 294,
Kluwer, 1990,
pp. 257--292.

\refitem{[24]}
T. H. Koornwinder,
{\sl The addition formula for little $q$-Legendre polynomials and the
$SU(2)$ quantum group},
SIAM J. Math. Anal. 22 (1991), 295--301.

\refitem{[25]}
T. H. Koornwinder,
{\sl Compact quantum groups and $q$-special functions},
in {\sl Representations of Lie groups and quantum groups},
V. Baldoni \& M. A. Picardello (eds.),
Pitman Research Notes in Mathematics Series 311,
Longman Scientific \& Technical, pp. 46--128.

\refitem{[26]}
S. Majid,
{\sl C-statistical quantum groups and Weyl algebras},
J. Math. Phys. 33 (1992), 3431--3444.

\refitem{[27]}
S. Majid,
{\sl Algebras and Hopf algebras in braided categories},
in {\sl Advances in Hopf algebras},
Marcel Dekker, 1994, pp. 55--105.

\refitem{[28]}
R. J. McDermott \& A. I. Solomon,
{\sl An analogue of the unitary displacement operator for the $q$-oscillator},
J. Phys. A 27 (1994), 2037--2043.

\refitem{[29]}
G. Polya \& G. L. Alexanderson,
{\sl Gaussian binomial coefficients},
Elem. Math. 26 (1970), 102--108.

\refitem{[30]}
L. J. Rogers,
{\sl On the expansion of some infinite products},
Proc. London Math. Soc. 24 (1893), 337--352.

\refitem{[31]}
G.-C. Rota, D. Kahaner, A. Odlyzko,
{\sl Finite operator calculus},
J. Math. Anal. Appl. 42 (1973), 684--760.

\refitem{[32]}
M. P. Sch\"utzenberger,
{\sl Une interpr\'etation de certaines solutions de l'\'equation
fonctionnelle: $F(x+y)=F(x)F(y)$},
C. R. Acad. Sci. Paris 236 (1953), 352--353.

\refitem{[33]}
M. E. Sweedler,
{\sl Hopf algebras},
Benjamin, 1969.

\refitem{[34]}
K. Ueno \& M. Nishizawa,
{\sl Quantum groups and zeta functions},
in {\sl Quantum groups: formalism and applications},
J. Lukierski e.a. (eds.),
PWN, Warsaw, 1995.

\refitem{[35]}
V. S. Varadarajan,
{\sl Lie groups, Lie algebras and their representations},
Prentice-Hall, 1974.

\refitem{[36]}
S. L. Woronowicz,
{\sl Operator equalities related to quantum $E(2)$ group},
Comm. Math. Phys. 144 (1992), 417--428.

\bye